# MATRIX ITERATIONS WITH VERTICAL SUPPORT RESTRICTIONS

## DIEGO A. MEJÍA

ABSTRACT. We use coherent systems of FS iterations on a power set, which can be seen as matrix iteration that allows restriction on arbitrary subsets of the vertical component, to prove general theorems about preservation of certain type of unbounded families on definable structures and of certain mad families (like those added by Hechler's poset for adding an a.d. family) regardless of the cofinality of their size. In particular, we define a class of posets called $\sigma$-*Frechet-linked* and show that they work well to preserve mad families, and unbounded families on $\omega^\omega$.

As applications of this method, we show that a large class of FS iterations can preserve the mad family added by Hechler's poset (regardless of the cofinality of its size), and the consistency of a constellation of Cichoń's diagram with 7 values where two of these values are singular.

## 1. INTRODUCTION

**Background.** In the framework of FS (finite support) iterations of ccc posets to prove consistency results with large continuum (that is, the size of the continuum $\mathfrak{c} = 2^{\aleph_0}$ larger than $\aleph_2$), very recently in [FFMM18] appeared the general notion of *coherent systems of FS iterations* that was used to construct a three-dimensional array of ccc posets to force that the cardinals in *Cichoń's diagram* are separated into 7 different values (see Figure 1). This is the first example of a 3D iteration that was used to prove a new consistency result. Moreover, the methods from [BF11] where used there to force, in addition, that the *almost disjointness number* $\mathfrak{a}$ is equal to the *bounding number* $\mathfrak{b}$, and to expand well-known results about preservation of mad families along FS iterations.

For quite some time, consistency results about many different values for cardinal invariants has been investigated. Some of the earliest results are due to Brendle [Bre91] who fixed standard techniques for FS iterations in this direction, and due to Blass and Shelah [BS89] who constructed the first example of a two-dimensional array of ccc posets to force the consistency of the existence of a base for a non-principal ultrafilter in $\omega$ of size smaller than the *dominating number* $\mathfrak{d}$. The latter technique received the name *matrix iterations* in [BF11] and it was improved there to prove the consistency of e.g. $\aleph_1 < \mathfrak{b} = \mathfrak{a} < \mathfrak{s}$. Recent developments on matrix iterations appear in work of the author [Mej13a, Mej13b], where forcing models satisfying that several cardinals in Cichoń's diagram are pairwise different (at most 6 different values were achieved) are constructed, and of Dow and Shelah [DS18] where the splitting number $\mathfrak{s}$ is forced to be singular. Concerning Cichoń's diagram, a few months ago Goldstern, Kellner and Shelah [GKS] used Boolean ultrapowers of strongly compact cardinals applied to the iteration constructed in [GMS16] to

2010 *Mathematics Subject Classification.* 03E17, 03E15, 03E35, 03E40.

Supported by grant no. IN201711, Dirección Operativa de Investigación, Institución Universitaria Pascual Bravo, and by the Grant-in-Aid for Early Career Scientists 18K13448, Japan Society for the Promotion of Science.

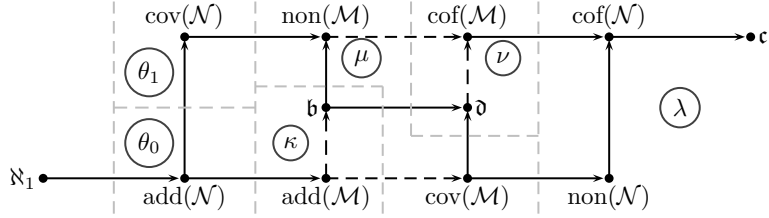

Figure 1. 7 values in Cichoń's diagram.

prove the consistency, modulo the existence of 4 strongly compact cardinals, of a division of Cichoń's diagram into 10 different values (the maximum number of values allowed in the diagram). Another example with 10 values, modulo 4 strongly compact cardinals, appears in [KST]. At this point, it is still unknown how to prove the consistency of 8 different values in this diagram modulo ZFC alone.

One drawback of the methods discussed so far is that the posets they produce force that the cardinal invariants that are not equal to $\mathfrak{c}$ must be regular. In the context of ccc forcing, one of the few exceptions is the consistency of $\mathfrak{d} < \mathfrak{c}$ with both cardinals singular, which can be obtained by a FS iteration where the last iterand is a random algebra (see e.g. [FFMM18, Thm. 5.1(d)]). On the other hand, many examples can be obtained by creature forcing constructions as in [KS09, KS12, FGKS17], for instance, in the latter reference it is proved that the right side of Cichoń's diagram can be divided into 5 different values where 4 of them are singular. However, all these constructions are $\omega^\omega$-bounding, so they force $\mathfrak{d} = \aleph_1$ and do not allow separation of cardinal invariants below $\mathfrak{d}$.

**Objective 1.** The main motivation of this research is to improve some of the ccc forcing methods to produce models where many cardinal invariants of the continuum are different and two or more of them are singular. As one of the main results of this paper, we show how to take advantage of the generality of coherent systems of FS iterations to produce such models where 2 cardinal invariants can be forced to be singular. In particular, we show that the 3D iteration of [FFMM18] that forces the constellation of 7 values in Cichoń's diagram can be modified so that 2 cardinals are allowed to be singular (Theorem 4.3). In addition. we modify examples from [Mej13a, FFMM18] in the same way.

**Methods.** A coherent system of FS iterations of length $\pi$ consists of a partial order $\langle I, \leq \rangle$ and, for each $i \in I$, a FS iteration $\mathbb{P}_{i,\pi} := \langle \mathbb{P}_{i,\alpha}, \dot{\mathbb{Q}}_{i,\alpha} : \alpha < \pi \rangle$ such that any pair of such iterations are coherent in the sense that, whenever $i \leq j$ in $I$ and $\alpha \leq \pi$, the $\mathbb{P}_{i,\alpha}$-generic extension is contained in the $\mathbb{P}_{j,\alpha}$-generic extension (see details in Definition 2.3 and a picture in Figure 6). For instance, a matrix iteration is a coherent system (of FS iterations) when $\langle I, \leq \rangle$ is a well-order (see Figure 2), and a 3D iteration is a coherent system on a product of ordinals $I = \gamma \times \delta$ with the coordinate-wise order (see Figure 3).

For our applications, we construct coherent systems on partial orders of the form $\langle \mathcal{P}(\Omega), \subseteq \rangle$, which in fact look like matrix iterations, with vertical component indexed by $\Omega$, that allow restriction on any arbitrary subset of $\Omega$. To be more precise, as the final generic extension of the forcing produced by such a system comes from the FS iteration $\langle \mathbb{P}_{\Omega,\xi}, \dot{\mathbb{Q}}_{\Omega,\xi} : \xi < \pi \rangle$, for any $A \subseteq \Omega$ the iteration $\langle \mathbb{P}_{A,\xi}, \dot{\mathbb{Q}}_{A,\xi} : \xi < \pi \rangle$ can be understood as the 'vertical' restriction on $A$ of the former FS iteration (see Figure 4). This "restriction" feature is what allows nice combinatorial arguments to force singular values for some cardinal invariants. Concretely, it allows to preserve *unbounded reals* (with respect to general structures, see Definition 3.4) that come from the vertical component (Theorem 3.15), and even maximal almost disjoint (mad) families of size of singular cardinality (Theorem 3.32). Surprisingly, the three-dimensional forcings from [FFMM18] can be reconstructed now as matrix iterations with vertical support restriction, though the real picture of the latter is the 'shape' of the lattice $\langle \mathcal{P}(\Omega), \subseteq \rangle$ plus one additional dimension (for the FS iterations).

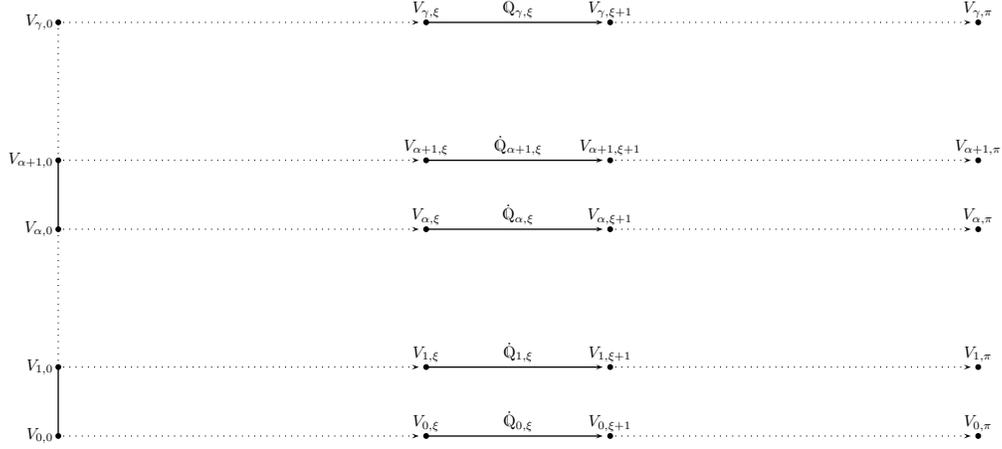

Figure 2. Matrix iteration.

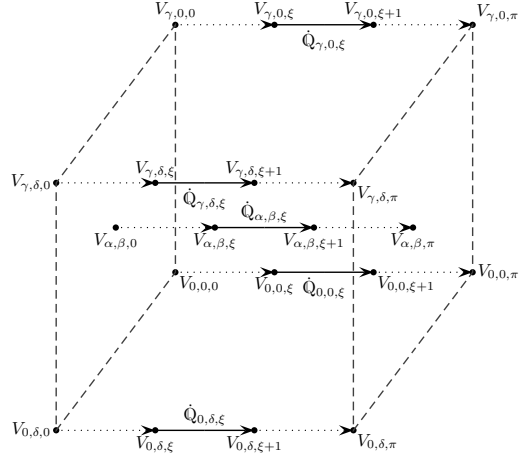

Figure 3. Three-dimensional iteration.

**Objective 2.** The theory of Brendle and Fischer [BF11] for preserving mad families is the cornerstone for the preservation results we propose in this paper, as well as it is in [FFMM18]. In the latter reference it is proved that $\mathbb{E}$ (the standard $\sigma$-centered poset adding an eventually different real) and random forcing (thus any random algebra) behaves well in their preservation theory, which allows to prove in [FFMM18, Thm. 4.17] that, whenever $\kappa$ is an uncountable regular cardinal, the mad family added by the Hechler poset $\mathbb{H}_\kappa$ is preserved by any further FS iteration whose iterands are either $\mathbb{E}$, a random algebra or a ccc poset of size $< \kappa$. In relation to this, we define a class of posets, which we call $\sigma$-*Frechet-linked* (see Definition 3.24), that includes $\mathbb{E}$ and random forcing, and we prove that any (definable) poset in this class behaves well with Brendle's and Fischer's preservation theory (Theorem 3.27). Moreover, by using coherent systems on a power set $\langle \mathcal{P}(\Omega), \subseteq \rangle$ we generalize [FFMM18, Thm. 4.17] by proving that, whenever $\Omega$ is uncountable (not necessarily of regular size), the mad family added by $\mathbb{H}_\Omega$ can be preserved by a large class of FS iterations (which includes the Suslin $\sigma$-Frechet-linked posets as iterands, see Theorem 4.1). This is related to the preservation of mad families of singular size discussed in the previous paragraph. In addition we also show that, for a cardinal $\mu$, $\mu$-Frechet-linked posets behave well in the preservation theory of unbounded families (Theorem 3.30).

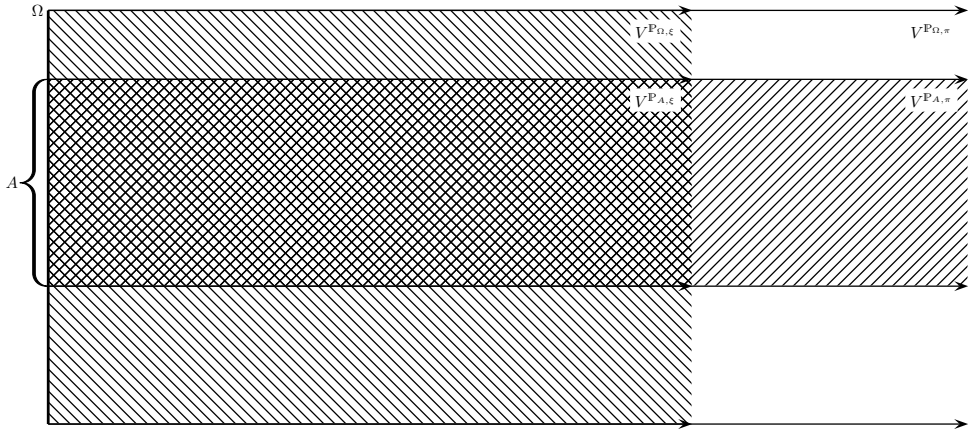

FIGURE 4. Matrix iteration with vertical support restriction.

**Structure of the paper.** In Section 2 we review the notion of coherent systems of FS iterations and prove general theorems about (vertical) direct limits within such a system. Section 3 is divided in two parts. In the first part, we review Judah's and Shelah's [JS90] and Brendle's [Bre91] theory of preservation of strongly unbounded families (with respect to general definable structures), as well as known facts from [BS89, BF11, Mej13a] to preserve unbounded reals. At the end, a general theorem about preservation of unbounded reals through coherent systems on a power set $\langle \mathcal{P}(\Omega), \subseteq \rangle$ is proved. In the second part we review Brendle's and Fischer's theory for mad family preservation, define $\mu$-Frechet-linked posets and prove that they behave well in this preservation theory. This allows to prove at the end a general theorem about preservation of mad families through coherent systems on a power set $\langle \mathcal{P}(\Omega), \subseteq \rangle$. Afterwards, in Section 4 we show applications of the theory presented so far, namely, mad family preservation along a large class of FS iterations and consistency results about Cichoń's diagram.

The last section proposes a general framework for linkedness of subsets of posets that includes notions like $n$-linked, centered and Frechet-linked. We say that $\Gamma$ is a *linkedness*

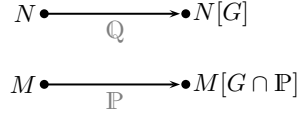

FIGURE 5. Generic extensions of pairs of posets ordered like $\mathbb{P} \lessdot_M \mathbb{Q}$.

*property of subsets of posets* if $\Gamma(\mathbb{P})$ defines a family of subsets of $\mathbb{P}$ for each poset $\mathbb{P}$. For such a property $\Gamma$ we define its corresponding notions of $\theta$-$\Gamma$-*Knaster* (a poset $\mathbb{P}$ has this property iff any subset of size $\theta$ contains a subset in $\Gamma(\mathbb{P})$ of size $\theta$) and $\mu$-$\Gamma$-*covered* (the version of $\mu$-linked for $\Gamma$). Built on the classical FS product and iteration theorems for Knaster and $\sigma$-linked, we find sufficient conditions for $\Gamma$ to generalize these theorems for $\theta$-$\Gamma$-Knaster and $\mu$-$\Gamma$-covered.

**Some notation.** Denote by $\mathbb{C}$ the poset that adds one Cohen real and by $\mathbb{C}_\Omega$ the poset that adds a family of Cohen reals indexed by the set $\Omega$ (which is basically a finite support product of $\mathbb{C}$). The Lebesgue measure on the Cantor space $2^\omega$ is denoted by Lb. *Random forcing*, denoted by $\mathbb{B}$, is the poset of Borel subsets of $2^\omega$ with positive Lebesgue measure, ordered by $\subseteq$. A *random algebra* on a set $\Omega$, denoted by $\mathbb{B}_\Omega$, is the poset of subsets of $2^{\Omega \times \omega}$ of the form $B \times 2^{(\Omega \smallsetminus J) \times \omega}$ for some $J \subseteq \Omega$ countable and some Borel subset $B$ of $2^{J \times \omega}$ with positive Lebesgue measure, ordered by $\subseteq$. This adds a family of random reals indexed by $\Omega$. Hechler poset for adding a dominating real is denoted by $\mathbb{D}$, and $\mathbb{E}$ is defined as the poset whose conditions are pairs $(s, \varphi)$ with $s \in \omega^{<\omega}$ and $\varphi : \omega \to [\omega]^{\leq m}$ for some $m < \omega$, ordered by $(s', \varphi') \leq (s, \varphi)$ iff $s \subseteq s'$, $\varphi(i) \subseteq \varphi'(i)$ for any $i < \omega$, and $s'(i) \notin \varphi(i)$ for any $i \in |s'| \smallsetminus |s|$. The trivial poset is denoted by $\mathbb{1}$.

Most of the cardinal invariants used in this paper are defined (or characterized) in Example 3.7. Recall that $A \subseteq [\omega]^{\aleph_0}$ is an *almost disjoint (a.d.) family* if the intersection of any two different members of $A$ is finite. A mad family is a maximal a.d. family, and $\mathfrak{a}$ is defined as the smallest size of an infinite mad family. For a set $\Omega$, Hechler's poset $\mathbb{H}_\Omega$ for adding an a.d. family (indexed by $\Omega$) is defined as the poset whose conditions are of the form $p : F_p \times n_p \to 2$ with $F_p \in [\Omega]^{<\omega}$ and $n_p < \omega$ (demand $n_p = 0$ iff $F_p = \emptyset$), ordered by $q \leq p$ iff $p \subseteq q$ and $|q^{-1}[\{1\}] \cap (F_p \times \{i\})| \leq 1$ for every $i \in [n_p, n_q)$ (see [Hec72]). This poset has the Knaster property and the a.d. family it adds is maximal when $\Omega$ is uncountable. It is forcing equivalent to $\mathbb{C}$ when $\Omega$ is countable and non-empty, and it is equivalent to $\mathbb{C}_{\omega_1}$ when $|\Omega| = \aleph_1$. For any $\Omega \subseteq \Omega'$, $\mathbb{H}_\Omega \lessdot \mathbb{H}_{\Omega'}$.

## 2. COHERENT SYSTEMS OF FS ITERATIONS

**Definition 2.1.** Let $M$ be a transitive model of ZFC. When $\mathbb{P} \in M$ and $\mathbb{Q}$ are posets, say that $\mathbb{P}$ *is a complete subposet of $\mathbb{Q}$ with respect to $M$*, abbreviated $\mathbb{P} \lessdot_M \mathbb{Q}$, if $\mathbb{P}$ is a subposet of $\mathbb{Q}$ and any maximal antichain of $\mathbb{P}$ that belongs to $M$ is still a maximal antichain in $\mathbb{Q}$.

If in addition $N$ is another transitive model of ZFC, $M \subseteq N$ and $\mathbb{Q} \in N$, then $\mathbb{P} \lessdot_M \mathbb{Q}$ implies that, whenever $G$ is $\mathbb{Q}$-generic over $N$, $G \cap \mathbb{P}$ is $\mathbb{P}$-generic over $M$ and $M[G \cap \mathbb{P}] \subseteq N[G]$ (see Figure 5).

**Example 2.2.** Let $M \subseteq N$ be transitive models of ZFC. When $\mathbb{P} \in M$ it is clear that $\mathbb{1} \lessdot_M \mathbb{P}$ and $\mathbb{P} \lessdot_M \mathbb{P}$. Also, if $\mathbb{S}$ is a Suslin ccc poset or a random algebra coded in $M$ then $\mathbb{S}^M \lessdot_M \mathbb{S}^N$.

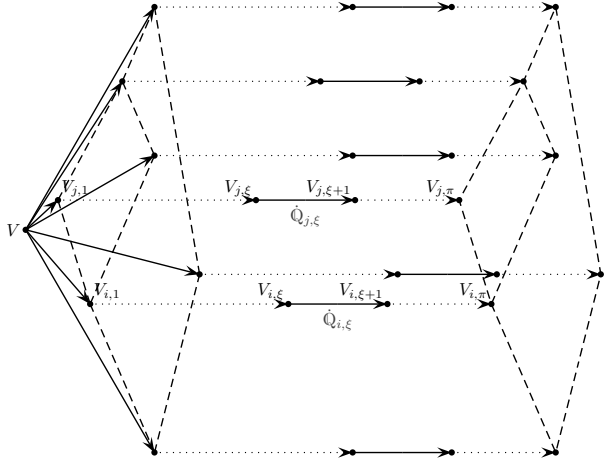

FIGURE 6. Coherent system of FS iterations. The figures in dashed lines represent the 'shape' of the partial order $\langle I, \leq \rangle$.

**Definition 2.3** ([FFMM18, Def. 3.2]). A *coherent system (of FS iterations)* $\mathbf{s}$ is composed of the following objects:

(I) a partially ordered set $I^{\mathbf{s}}$, an ordinal $\pi^{\mathbf{s}}$, and

(II) for each $i \in I^{\mathbf{s}}$, a FS iteration $\mathbb{P}^{\mathbf{s}}_{i, \pi^{\mathbf{s}}} = \langle \mathbb{P}^{\mathbf{s}}_{i,\xi}, \dot{\mathbb{Q}}^{\mathbf{s}}_{i,\xi} : \xi < \pi^{\mathbf{s}} \rangle$ such that, for any $i \leq j$ in $I$ and $\xi < \pi^{\mathbf{s}}$, if $\mathbb{P}^{\mathbf{s}}_{i,\xi} \lessdot \mathbb{P}^{\mathbf{s}}_{j,\xi}$ then $\mathbb{P}^{\mathbf{s}}_{j,\xi}$ forces $\dot{\mathbb{Q}}^{\mathbf{s}}_{i,\xi} \lessdot_{V^{\mathbb{P}^{\mathbf{s}}_{i,\xi}}} \dot{\mathbb{Q}}^{\mathbf{s}}_{j,\xi}$.

According to this notation, $\mathbb{P}^{\mathbf{s}}_{i,0}$ is the trivial poset and $\mathbb{P}^{\mathbf{s}}_{i,1} = \dot{\mathbb{Q}}^{\mathbf{s}}_{i,0}$. We often refer to $\langle \mathbb{P}^{\mathbf{s}}_{i,1} : i \in I^{\mathbf{s}} \rangle$ as the *base of the coherent system* $\mathbf{s}$. The condition given in (II) implies that $\mathbb{P}^{\mathbf{s}}_{i,\xi} \lessdot \mathbb{P}^{\mathbf{s}}_{j,\xi}$ whenever $i \leq j$ in $I^{\mathbf{s}}$ and $\xi \leq \pi^{\mathbf{s}}$ (see Lemma 3.14).

For $j \in I^{\mathbf{s}}$ and $\eta \leq \pi^{\mathbf{s}}$ we write $V^{\mathbf{s}}_{j,\eta}$ for the $\mathbb{P}^{\mathbf{s}}_{j,\eta}$-generic extensions. Concretely, when $G$ is $\mathbb{P}^{\mathbf{s}}_{j,\eta}$-generic over $V$, $V^{\mathbf{s}}_{j,\eta} := V[G]$ and $V^{\mathbf{s}}_{i,\xi} := V[\mathbb{P}^{\mathbf{s}}_{i,\xi} \cap G]$ for all $i \leq j$ in $I^{\mathbf{s}}$ and $\xi \leq \eta$. Note that $V^{\mathbf{s}}_{i,\xi} \subseteq V^{\mathbf{s}}_{j,\eta}$ and $V^{\mathbf{s}}_{i,0} = V$ (see Figure 6).

We say that the coherent system $\mathbf{s}$ has the *ccc* if, additionally, $\mathbb{P}^{\mathbf{s}}_{i,\xi}$ forces that $\dot{\mathbb{Q}}^{\mathbf{s}}_{i,\xi}$ has the ccc for each $i \in I^{\mathbf{s}}$ and $\xi < \pi^{\mathbf{s}}$. This implies that $\mathbb{P}^{\mathbf{s}}_{i,\xi}$ has the ccc for all $i \in I^{\mathbf{s}}$ and $\xi \leq \pi^{\mathbf{s}}$.

A concrete simple type of coherent system is what we call a *coherent pair (of FS iterations)*. A coherent system $\mathbf{s}$ is a coherent pair if $I^{\mathbf{s}}$ is of the form $\{i_0, i_1\}$ ordered as $i_0 < i_1$.

For a coherent system $\mathbf{s}$ and a set $J \subseteq I^{\mathbf{s}}$, $\mathbf{s}|J$ denotes the coherent system with $I^{\mathbf{s}|J} = J$, $\pi^{\mathbf{s}|J} = \pi^{\mathbf{s}}$ and the FS iterations corresponding to (II) defined as for $\mathbf{s}$; if $\eta \leq \pi^{\mathbf{s}}$, $\mathbf{s}|\eta$ denotes the coherent system with $I^{\mathbf{s}|\eta} = I^{\mathbf{s}}$, $\pi^{\mathbf{s}|\eta} = \eta$ and the iterations for (II) defined up to $\eta$ as for $\mathbf{s}$. Note that, if $i_0 < i_1$ in $I^{\mathbf{s}}$, then $\mathbf{s}|\{i_0, i_1\}$ is a coherent pair and $\mathbf{s}|\{i_0\}$ is just the FS iteration $\mathbb{P}^{\mathbf{s}}_{i,\pi^{\mathbf{s}}} = \langle \mathbb{P}^{\mathbf{s}}_{i,\xi}, \dot{\mathbb{Q}}^{\mathbf{s}}_{i,\xi} : \xi < \pi^{\mathbf{s}} \rangle$.

In particular, the upper indices $\mathbf{s}$ are omitted when there is no risk of ambiguity.

The following is a generalization of [FFMM18, Lemma 3.7].

**Lemma 2.4.** *Let $\theta$ be an uncountable regular cardinal. Assume that $\mathbf{s}$ is a coherent system that satisfies:*

(i) $I$ has a maximum $i^*$ and $I \smallsetminus \{i^*\}$ is $< \theta$-directed,

(ii) each $\mathbb{P}_{i,\xi}$ forces that $\dot{\mathbb{Q}}_{i,\xi}$ is $\theta$-cc, and

(iii) for any $\xi < \pi$, if $\mathbb{P}_{i^*,\xi}$ is the direct limit of $\langle \mathbb{P}_{i,\xi} : i < i^* \rangle$ then $\mathbb{P}_{i^*,\xi}$ forces that $\dot{\mathbb{Q}}_{i^*,\xi} = \bigcup_{i<i^*} \dot{\mathbb{Q}}_{i,\xi}$.

Then, for any $\xi \leq \pi$,

(a) $\mathbb{P}_{i^*,\xi}$ is the direct limit of $\langle \mathbb{P}_{i,\xi} : i < i^* \rangle$ and

(b) if $\gamma < \theta$ and $\dot{f}$ is a $\mathbb{P}_{i^*,\xi}$-name of a function from $\gamma$ into $\bigcup_{i<i^*} V_{i,\xi}$ then $\dot{f}$ is (forced to be equal to) a $\mathbb{P}_{i,\xi}$-name for some $i < i^*$. In particular, the reals in $V_{i^*,\xi}$ are precisely the reals in $\bigcup_{i<i^*} V_{i,\xi}$.

*Proof.* By condition (ii) it is clear that $\mathbb{P}_{i^*,\pi}$ has the $\theta$-cc. We first show that (a) implies (b). Let $\dot{f}$ be as in the condition in (b). For each $\alpha < \gamma$ choose a maximal antichain $A_\alpha$ that decides some $i < i^*$ and some $\mathbb{P}_{i,\xi}$-name to which $\dot{f}(\alpha)$ is forced to be equal. As $|A_\alpha| < \theta$ and $I \smallsetminus \{i^*\}$ is $< \theta$-directed, by (a) we can find some $i_\alpha < i^*$ and some $\mathbb{P}_{i_\alpha,\xi}$-name $\dot{x}_\alpha$ such that $\dot{f}(\alpha)$ is forced to be equal to $\dot{x}_\alpha$. Hence, as $\gamma < \theta$, there is some upper bound $j < i^*$ of $\{i_\alpha : \alpha < \gamma\}$, so $\dot{f}$ is forced to be equal to the $\mathbb{P}_{j,\xi}$-name of the function $\alpha \mapsto \dot{x}_\alpha$.

Now we prove (a) by induction on $\xi$. The case $\xi = 0$ and the limit step are clear. Assume that (a) holds for $\xi$, and we show that (a) holds for $\xi + 1$. If $p \in \mathbb{P}_{i^*,\xi+1}$ then, by (i), (iii) and (b), there is some $j_1 < i^*$ such that $p(\xi)$ is (forced to be equal to) a $\mathbb{P}_{j_1,\xi}$-name of a member of $\dot{\mathbb{Q}}_{j_1,\xi}$. On the other hand, by (a), there is some $j_0 < i^*$ such that $p{\restriction}\xi \in \mathbb{P}_{j_0,\xi}$. Hence, by (i), there is some $j < i^*$ above $j_0$ and $j_1$, so $p \in \mathbb{P}_{j,\xi}$. $\qquad\square$

Typically, a coherent system of FS iterations is constructed by transfinite recursion. Concretely, to construct such a system $\mathbf{s}$ of FS iterations of length $\pi$, $\mathbf{s}{\restriction}\xi$ is constructed by recursion on $\xi \leq \pi$ as follows. In the step $\xi = 0$, we determine the partial order $\langle I, \leq \rangle$ that will support the base of the coherent system; in the limit step it is just enough to take direct limits of the systems constructed previously; for the successor step, assuming that $\mathbf{s}{\restriction}\xi$ has been constructed, the system $\langle \dot{\mathbb{Q}}_{i,\xi} : i \in I \rangle$ of names of posets ($\dot{\mathbb{Q}}_{i,\xi}$ is a $\mathbb{P}_{i,\xi}$-name) that will determine how it is forced in stage $\xi$ is determined, and afterwards $\mathbf{s}{\restriction}(\xi+1)$ is defined so that it extends $\mathbf{s}{\restriction}\xi$ and $\mathbb{P}_{i,\xi+1} = \mathbb{P}_{i,\xi} * \dot{\mathbb{Q}}_{i,\xi}$ for each $i \in I$. To have that $\mathbf{s}{\restriction}(\xi+1)$ is indeed a coherent system, we require that $\Vdash_{\mathbb{P}_{i,\xi}} \dot{\mathbb{Q}}_{i,\xi} \lessdot_{V_{i,\xi}} \dot{\mathbb{Q}}_{j,\xi}$ whenever $i \leq j$ in $I$.

In short, to construct a coherent system as in the previous paragraph, it is just enough to determine the partial order that will serve as base and to determine the iterands suitably. As an example (that will serve for all our applications), we define the following simple type of coherent systems. Suslin posets play an important role in such systems (see [JS88]).

**Definition 2.5** ([FFMM18, Def. 3.8] with variations). A coherent system $\mathbf{s}$ is *standard* if

(I) it consists, additionally, of

(i) a partition $\langle S^{\mathbf{s}}, C^{\mathbf{s}} \rangle$ of $\pi^{\mathbf{s}} \smallsetminus \{0\}$,

(ii) a function $\Delta^{\mathbf{s}} : [1, \pi^{\mathbf{s}}) \to I^{\mathbf{s}}$,

(iii) a sequence $\langle \dot{S}^{\mathbf{s}}_\xi : \xi \in S^{\mathbf{s}} \rangle$ where each $\dot{S}^{\mathbf{s}}_\xi$ is either a (name of a definition of a) Suslin ccc poset coded in $V^{\mathbf{s}}_{\Delta(\xi),\xi}$, or a random algebra, and

(iv) a sequence $\langle \dot{\mathbb{Q}}^{\mathbf{s}}_\xi : \xi \in C^{\mathbf{s}} \rangle$ such that each $\dot{\mathbb{Q}}^{\mathbf{s}}_\xi$ is a $\mathbb{P}^{\mathbf{s}}_{\Delta^{\mathbf{s}}(\xi),\xi}$-name of a poset that is forced to have the ccc by $\mathbb{P}^{\mathbf{s}}_{i,\xi}$ for every $i \geq \Delta^{\mathbf{s}}(\xi)$ in $I^{\mathbf{s}}$, and

(II) for any $i \in I^{\mathbf{s}}$ it satisfies:
  (i) $\dot{\mathbb{Q}}^{\mathbf{s}}_{i,0}$ has the ccc and
  (ii) for any $0 < \xi < \pi^{\mathbf{s}}$,

$$\dot{\mathbb{Q}}^{\mathbf{s}}_{i,\xi} = \begin{cases} (\dot{\mathbb{S}}^{\mathbf{s}}_\xi)^{V^{\mathbf{s}}_{i,\xi}} & \text{if } \xi \in S^{\mathbf{s}} \text{ and } i \geq \Delta^{\mathbf{s}}(\xi), \\ \dot{\mathbb{Q}}^{\mathbf{s}}_\xi & \text{if } \xi \in C^{\mathbf{s}} \text{ and } i \geq \Delta^{\mathbf{s}}(\xi), \\ \mathbb{1} & \text{otherwise.} \end{cases}$$

As in Definition 2.3, the upper index $\mathbf{s}$ may be omitted when understood.

Note that any standard coherent system has the ccc. One of the main features of a standard coherent system is the type of generic objects that are added, which is determined in principle by the partition $\langle S^{\mathbf{s}}, C^{\mathbf{s}} \rangle$ of $\pi^{\mathbf{s}} \smallsetminus \{0\}$. When $\xi \in S^{\mathbf{s}}$, the generic object added at stage $\xi + 1$ is called *full generic*, while in the case $\xi \in C^{\mathbf{s}}$, such generic object is called *restricted generic*. The reason for this is that, in the case $\xi \in S^{\mathbf{s}}$ the generic object added by $\dot{\mathbb{S}}_\xi$ is generic over $V_{i,\xi}$ *for all* $i \in I$, while in the case $\xi \in C^{\mathbf{s}}$ the generic object added by $\dot{\mathbb{Q}}_\xi$ is only generic over $V_{\Delta(\xi),\xi}$. For instance, this distinction is fundamental to deal with forcing constructions to separate several cardinal invariants of the continuum (as in our applications). Even more, a restricted generic can be even much more restricted, that is, not really generic over $V_{\Delta(\xi),\xi}$ but over some ZFC model $N \subseteq V_{\Delta(\xi),\xi}$ (e.g. when $\dot{\mathbb{Q}}_\xi = \mathbb{D}^N$).

The version of Lemma 2.4 for standard coherent systems requires simpler conditions.

**Corollary 2.6.** *Let $\mathbf{s}$ be a standard coherent system and let $\theta$ be an uncountable regular cardinal. If*

  (i) *$I$ has a maximum $i^*$, $I \smallsetminus \{i^*\}$ is $< \theta$-directed,*
  (ii) *$i^* \notin \mathrm{ran}\Delta$, and*
  (iii) *whenever $\pi > 0$, $\mathbb{P}_{i^*,1}$ is the direct limit of $\langle \mathbb{P}_{i,1} : i < i^* \rangle$,*

*then (a) and (b) of Lemma 2.4 hold.*

*Proof.* It is clear that hypotheses (i)-(iii) of Lemma 2.4 are satisfied. $\square$

In our applications we use coherent systems on a power set $\langle \mathcal{P}(\Omega), \subseteq \rangle$. If $\mathbf{s}$ is a standard coherent system based in such a partial order, then $\mathbb{P}_{\Omega,\pi}$ (the largest poset in the system) can be represented as a two dimensional forcing construction supported in the plane $\Omega \times \pi$ and, for any $A \subseteq \gamma$ and $\xi \leq \pi$, $\mathbb{P}_{A,\xi}$ is seen as the restriction of the construction to the rectangle $A \times \xi$ (see Figure 4 in Section 1). The following result is a suitable consequence of Corollary 2.6 when dealing with such (standard) coherent systems.

**Lemma 2.7.** *Let $\theta$ be a cardinal of uncountable cofinality and let $\mathbf{s}$ be a standard coherent system where $I = I^{\mathbf{s}}$ is a suborder of $\langle \mathcal{P}(\Omega), \subseteq \rangle$. Assume that*

  (i) *$I$ is closed under intersections,*
  (ii) *$I \cap [\Omega]^{<\theta}$ is cofinal in $[\Omega]^{<\theta}$,*
  (iii) *$\Delta(\xi) \in [\Omega]^{<\theta}$ for any $\xi \in [1, \pi)$ (see Definition 2.5(I)(ii)), and*
  (iv) *whenever $\pi > 0$ and $X \in I$, $\mathbb{P}_{X,1}$ is the direct limit of $\langle \mathbb{P}_{A,1} : A \in I \cap [X]^{<\theta} \rangle$.*

*Then, for every $X \in I$ and $\xi \leq \pi$,*

(a) *$\mathbb{P}_{X,\xi}$ is the direct limit of $\langle \mathbb{P}_{A,\xi} : A \in I \cap [X]^{<\theta} \rangle$ and*
(b) *for any $\mathbb{P}_{X,\xi}$-name of a function $\dot{x}$ with domain $\gamma < \mathrm{cf}(\theta)$ into $\bigcup_{A \in I \cap [X]^{<\theta}} V_{A,\xi}$, there is some $A \in I \cap [X]^{<\theta}$ such that $\dot{x}$ is (forced to be equal to) a $\mathbb{P}_{A,\xi}$-name.*

*Proof.* Fix $X \in I$. The lemma is trivial when $|X| < \theta$, so assume that $|X| \geq \theta$. By (i) and (ii), $I \cap [X]^{<\theta}$ is cofinal in $[X]^{<\theta}$. Put $I^* = I \cap ([X]^{<\theta} \cup \{X\})$. Hence $X$ is the maximum of $I^*$ and $I^* \smallsetminus \{X\} = I \cap [X]^{<\theta}$ is $< \mathrm{cf}(\theta)$-directed by (i) and (ii). Consider $\mathbf{s}^* = \mathbf{s}|I^*$. Note that $\mathbf{s}^*$ is a standard coherent system similar to $\mathbf{s}$ with the difference that $\Delta^* = \Delta^{\mathbf{s}^*} : [1, \pi) \to I \cap [X]^{<\theta}$ is defined as $\Delta^*(\xi) := \Delta(\xi)$ whenever $\Delta(\xi) \subseteq X$, or $\Delta^*(\xi) := \emptyset$ otherwise. Also, $\dot{\mathbb{S}}_\xi^{\mathbf{s}^*} = \dot{\mathbb{S}}_\xi^{\mathbf{s}}$ and $\dot{\mathbb{Q}}_\xi^{\mathbf{s}^*} = \dot{\mathbb{Q}}_\xi^{\mathbf{s}}$ in the first case, otherwise each one is the trivial poset. As $X \notin \mathrm{ran}\Delta^*$, the result is a direct consequence of Corollary 2.6 applied to $\mathbf{s}^*$ and $\mathrm{cf}(\theta)$. $\square$

**Example 2.8.** Let $\Omega$ be a set.

(1) The partial order $\langle \mathcal{P}(\Omega), \subseteq \rangle$ clearly satisfies conditions (i) and (ii) of Lemma 2.7 (for any infinite cardinal $\theta$).

Assume that $\Omega = \Omega_0 \cup \Omega_1$ is a disjoint union. If $\mathbf{s}$ is a standard coherent system on $\langle \mathcal{P}(\Omega), \subseteq \rangle$ such that $\mathbb{P}_{X,1} = \mathbb{H}_{X \cap \Omega_0} \times \mathbb{C}_{X \cap \Omega_1}$ for any $X \subseteq \Omega$, then condition (iv) of Lemma 2.7 is satisfied for $\theta = \aleph_1$ (and hence for any uncountable $\theta$).

(2) If $\theta$ is a regular cardinal, $\Omega = \theta$ and $I_0$ is a cofinal subset of $\theta$, then $I := I_0 \cup \{\theta\}$ satisfies conditions (i) and (ii) of Lemma 2.7. Such a partial order $I$ (in particular $I = \theta \cup \{\theta\}$) is used to construct classical matrix iterations as in, e.g., [BS89, BF11, Mej13a, FFMM18].

## 3. Preservation properties

As mentioned in the introduction, this section is divided in two parts. For convenience with the notation fixed in the first part, we use a different notation from [BF11, FFMM18] for the results in the second part.

3.1. **Preservation theory.** A generalization of the contents of this part, as well as complete proofs and more examples, can be found in [CM, Sect. 4].

Typically, cardinal invariants of the continuum are defined through *relational systems* as follows.

**Definition 3.1.** A *relational system* is a triplet $\mathbf{R} = \langle X, Y, \sqsubset \rangle$ where $\sqsubset$ is a relation contained in $X \times Y$. For $x \in X$ and $y \in Y$, $x \sqsubset y$ is often read $y \sqsubset$-*dominates* $x$. A family $F \subseteq X$ is $\mathbf{R}$-*unbounded* if there is <u>no</u> real in $Y$ that $\sqsubset$-dominates every member of $F$. Dually, $D \subseteq Y$ is an $\mathbf{R}$-*dominating* family if every member of $X$ is $\sqsubset$-dominated by some member of $D$. The cardinal $\mathfrak{b}(\mathbf{R})$ denotes the least size of an $\mathbf{R}$-unbounded family and $\mathfrak{d}(\mathbf{R})$ is the least size of an $\mathbf{R}$-dominating family.

Say that $x \in X$ is $\mathbf{R}$-*unbounded over a set* $M$ if $x \not\sqsubset y$ for all $y \in Y \cap M$. Given a cardinal $\lambda$ say that $F \subseteq X$ is $\lambda$-$\mathbf{R}$-*unbounded* if, for any $Z \subseteq Y$ of size $< \lambda$, there is an $x \in F$ that is $\mathbf{R}$-unbounded over $Z$. Say that $F \subseteq X$ is *strongly* $\lambda$-$\mathbf{R}$-*unbounded* if $|F| \geq \lambda$ and $|\{x \in F : x \sqsubset y\}| < \lambda$ for any $y \in Y$.

**Remark 3.2.** When $\lambda \geq 2$, any $\lambda$-$\mathbf{R}$-unbounded family is $\mathbf{R}$-unbounded. Hence, if $F$ is a $\lambda$-$\mathbf{R}$-unbounded family then $\mathfrak{b}(\mathbf{R}) \leq |F|$ and $\lambda \leq \mathfrak{d}(\mathbf{R})$. Also, if $\theta$ is regular and $F'$ is a strongly $\theta$-$\mathbf{R}$-unbounded family then it is $|F'|$-$\mathbf{R}$-unbounded, so $\mathfrak{b}(\mathbf{R}) \leq |F'| \leq \mathfrak{d}(\mathbf{R})$.

**Definition 3.3.** Let $\mathbf{R} = \langle X, Y, \sqsubset \rangle$ and $\mathbf{R}' = \langle X', Y', \sqsubset' \rangle$ be two relational systems. Say that $\mathbf{R}$ is *Tukey-Galois below* $\mathbf{R}'$ if there are two maps $F : X \to X'$ and $G : Y' \to Y$ such that, for each $x \in X$ and $b \in Y'$, if $F(x) \sqsubset' b$ then $x \sqsubset G(b)$. When, in addition, $\mathbf{R}'$ is Tukey-Galois below $\mathbf{R}$, we say that $\mathbf{R}$ and $\mathbf{R}'$ are *Tukey-Galois equivalent*.

Recall that, whenever $\mathbf{R}$ is Tukey-Galois below $\mathbf{R}'$, $\mathfrak{b}(\mathbf{R}') \leq \mathfrak{b}(\mathbf{R})$ and $\mathfrak{d}(\mathbf{R}) \leq \mathfrak{d}(\mathbf{R}')$.

**Definition 3.4.** A relational system $\mathbf{R} := \langle X, Y, \sqsubset \rangle$ is a *Polish relational system (Prs)* if the following is satisfied:

(i) $X$ is a perfect Polish space,
(ii) $Y$ is a non-empty analytic subspace of some Polish space $Z$ and
(iii) $\sqsubset = \bigcup_{n<\omega} \sqsubset_n$ for some increasing sequence $\langle \sqsubset_n \rangle_{n<\omega}$ of closed subsets of $X \times Z$ such that $(\sqsubset_n)^y = \{x \in X : x \sqsubset_n y\}$ is nwd (nowhere dense) for all $y \in Y$.

By (iii), $\langle X, \mathcal{M}(X), \in \rangle$ is Tukey-Galois below $\mathbf{R}$ where $\mathcal{M}(X)$ denotes the $\sigma$-ideal of meager subsets of $X$. Therefore, $\mathfrak{b}(\mathbf{R}) \leq \mathrm{non}(\mathcal{M})$ and $\mathrm{cov}(\mathcal{M}) \leq \mathfrak{d}(\mathbf{R})$. Moreover, (iii) implies that, whenever $c \in X$ is a Cohen real over a transitive model $M$ of ZFC and the Prs $\mathbf{R}$ is coded in $M$, $c$ is $\mathbf{R}$-unbounded over $M$.

**Definition 3.5** (Judah and Shelah [JS90])**.** Let $\mathbf{R} = \langle X, Y, \sqsubset \rangle$ be a Prs and let $\theta$ be a cardinal. A poset $\mathbb{P}$ is $\theta$-$\mathbf{R}$-*good* if, for any $\mathbb{P}$-name $\dot{h}$ for a real in $Y$, there is a non-empty $H \subseteq Y$ of size $< \theta$ such that $\Vdash x \not\sqsubset \dot{h}$ for any $x \in X$ (in the ground model) that is $\mathbf{R}$-unbounded over $H$.

Say that $\mathbb{P}$ is $\mathbf{R}$-*good* when it is $\aleph_1$-$\mathbf{R}$-good.

Definition 3.5 describes a property, respected by FS iterations, to preserve specific types of $\mathbf{R}$-unbounded families. Concretely, when $\theta$ is uncountable regular,

(a) any $\theta$-$\mathbf{R}$-good poset preserves all the (strongly) $\theta$-$\mathbf{R}$-unbounded families from the ground model and
(b) FS iterations of $\theta$-cc $\theta$-$\mathbf{R}$-good posets produce $\theta$-$\mathbf{R}$-good posets.

By Remark 3.2, posets that are $\theta$-$\mathbf{R}$-good work to preserve $\mathfrak{b}(\mathbf{R})$ small and $\mathfrak{d}(\mathbf{R})$ large.

Clearly, $\theta$-$\mathbf{R}$-good implies $\theta'$-$\mathbf{R}$-good whenever $\theta \leq \theta'$, and any poset completely embedded into a $\theta$-$\mathbf{R}$-good poset is also $\theta$-$\mathbf{R}$-good.

Consider the following particular cases of interest for our applications.

**Lemma 3.6** ([Mej13a, Lemma 4])**.** *If $\mathbf{R}$ is a Prs and $\theta$ is an uncountable regular cardinal, then any poset of size $< \theta$ is $\theta$-$\mathbf{R}$-good. In particular, Cohen forcing is $\mathbf{R}$-good.*

**Example 3.7.** Fix an uncountable regular cardinal $\theta$.

(1) *Preserving non-meager sets:* Consider the Polish relational system $\mathbf{Ed} := \langle \omega^\omega, \omega^\omega, \neq^* \rangle$ where $x \neq^* y$ iff $x$ and $y$ are eventually different, that is, $x(i) \neq y(i)$ for all but finitely many $i < \omega$. By [BJ95, Thm. 2.4.1 and 2.4.7], $\mathfrak{b}(\mathbf{Ed}) = \mathrm{non}(\mathcal{M})$ and $\mathfrak{d}(\mathbf{Ed}) = \mathrm{cov}(\mathcal{M})$.

(2) *Preserving unbounded families:* Let $\mathbf{D} := \langle \omega^\omega, \omega^\omega, \leq^* \rangle$ be the Polish relational system where $x \leq^* y$ iff $x(i) \leq y(i)$ for all but finitely many $i < \omega$. Clearly, $\mathfrak{b}(\mathbf{D}) = \mathfrak{b}$ and $\mathfrak{d}(\mathbf{D}) = \mathfrak{d}$.

  Miller [Mil81] proved that $\mathbb{E}$ is $\mathbf{D}$-good. Furthermore, $\omega^\omega$-bounding posets, like the random algebra, are $\mathbf{D}$-good. In Theorem 3.30 we prove that $\mu$-Frechet-linked posets are $\mu^+$-$\mathbf{D}$-good.

(3) *Preserving null-covering families:* Define $X_n := \{a \in [2^{<\omega}]^{<\aleph_0} : \mathrm{Lb}(\bigcup_{s \in a}[s]) \leq 2^{-n}\}$ (endowed with the discrete topology) and put $X := \prod_{n<\omega} X_n$ with the product topology, which is a perfect Polish space. For every $x \in X$ denote $N_x^* := \bigcap_{n<\omega} \bigcup_{s \in x_n}[s]$, which is clearly a Borel null set in $2^\omega$.

  Define the Prs $\mathbf{Cn} := \langle X, 2^\omega, \sqsubset \rangle$ where $x \sqsubset z$ iff $z \notin N_x^*$. Recall that any null set in $2^\omega$ is a subset of $N_x^*$ for some $x \in X$, so $\mathbf{Cn}$ and $\langle \mathcal{N}(2^\omega), 2^\omega, \not\ni \rangle$ are Tukey-Galois equivalent. Therefore $\mathfrak{b}(\mathbf{Cn}) = \mathrm{cov}(\mathcal{N})$ and $\mathfrak{d}(\mathbf{Cn}) = \mathrm{non}(\mathcal{N})$.

  By a similar argument as in [Bre91, Lemma 1*], any $\nu$-centered poset is $\theta$-$\mathbf{Cn}$-good for any $\nu < \theta$ infinite. In particular, $\sigma$-centered posets are $\mathbf{Cn}$-good.

(4) *Preserving "union of null sets is not null":* For each $k < \omega$ let $\mathrm{id}^k : \omega \to \omega$ such that $\mathrm{id}^k(i) = i^k$ for all $i < \omega$ and put $\mathcal{H} := \{\mathrm{id}^{k+1} : k < \omega\}$. Let $\mathbf{Lc} := \langle \omega^\omega, \mathcal{S}(\omega, \mathcal{H}), \in^* \rangle$ be the Polish relational system where

$$\mathcal{S}(\omega, \mathcal{H}) := \{\varphi : \omega \to [\omega]^{<\aleph_0} : \exists h \in \mathcal{H} \forall i < \omega(|\varphi(i)| \le h(i))\},$$

and $x \in^* \varphi$ iff $\exists n < \omega \forall i \ge n(x(i) \in \varphi(i))$, which is read *$x$ is localized by $\varphi$*. As a consequence of Bartoszyński's characterization (see [BJ95, Thm. 2.3.9]), $\mathfrak{b}(\mathbf{Lc}) = \mathrm{add}(\mathcal{N})$ and $\mathfrak{d}(\mathbf{Lc}) = \mathrm{cof}(\mathcal{N})$.

Any $\nu$-centered poset is $\theta$-$\mathbf{Lc}$-good for any $\nu < \theta$ infinite (see [JS90]) so, in particular, $\sigma$-centered posets are $\mathbf{Lc}$-good. Moreover, Kamburelis [Kam89] proved that any Boolean algebra with a strictly positive finitely additive measure is $\mathbf{Lc}$-good. As a consequence, subalgebras (not necessarily complete) of random forcing are $\mathbf{Lc}$-good.

**Lemma 3.8.** *If $\mathbf{R} = \langle X, Y, \sqsubset \rangle$ is a Prs then, for any set $\Omega$, $\mathbb{H}_\Omega$ is $\mathbf{R}$-good.*

*Proof.* Let $\dot{y}$ be a $\mathbb{H}_\Omega$-name of a member of $Y$. Then there is some countable $A \subseteq \Omega$ such that $\dot{y}$ is a $\mathbb{H}_A$-name. As $\mathbb{H}_A$ is countable, it is $\mathbf{R}$-good by Lemma 3.6, so there is some non-empty countable $H \subseteq Y$ witnessing this. The same $H$ witnesses goodness for $\dot{y}$ and $\mathbb{H}_\Omega$. □

In a similar way, it can be proved that any random algebra is $\mathbf{R}$-good iff random forcing is $\mathbf{R}$-good.

The following results indicate that (strongly) $\nu$-unbounded families can be added with Cohen reals, and the effect on $\mathfrak{b}(\mathbf{R})$ and $\mathfrak{d}(\mathbf{R})$ by a FS iteration of good posets.

**Lemma 3.9.** *Let $\nu$ be a cardinal of uncountable cofinality, $\mathbf{R} = \langle X, Y, \sqsubset \rangle$ a Prs and let $\langle \mathbb{P}_\alpha \rangle_{\alpha < \nu}$ be a $\lessdot$-increasing sequence of $\mathrm{cf}(\nu)$-cc posets such that $\mathbb{P}_\nu = \mathrm{limdir}_{\alpha < \nu} \mathbb{P}_\alpha$. If $\mathbb{P}_{\alpha+1}$ adds a Cohen real $\dot{c}_\alpha \in X$ over $V^{\mathbb{P}_\alpha}$ for any $\alpha < \nu$, then $\mathbb{P}_\nu$ forces that $\{\dot{c}_\alpha : \alpha < \nu\}$ is a strongly $\nu$-$\mathbf{R}$-unbounded family of size $\nu$.*

**Theorem 3.10.** *Let $\theta$ be an uncountable regular cardinal, $\mathbf{R} = \langle X, Y, \sqsubset \rangle$ a Prs, $\delta \ge \theta$ an ordinal and let $\langle \mathbb{P}_\alpha, \dot{\mathbb{Q}}_\alpha \rangle_{\alpha < \delta}$ be a FS iteration of non-trivial $\theta$-$\mathbf{R}$-good $\theta$-cc posets. Then, $\mathbb{P}_\delta$ forces $\mathfrak{b}(\mathbf{R}) \le \theta$ and $\mathfrak{d}(\mathbf{R}) \ge |\delta|$.*

*Proof.* See e.g. [CM, Thm. 4.15] or [GMS16, Cor. 3.6]. □

Fix transitive models $M \subseteq N$ of ZFC and a Polish relational system $\mathbf{R} = \langle X, Y, \sqsubset \rangle$ coded in $M$. The following results are related to preservation of $\mathbf{R}$-unbounded reals along coherent pairs of FS iterations.

**Lemma 3.11** ([Mej13a, Thm. 7]). *Let $\mathbb{S}$ be a Suslin ccc poset coded in $M$. If $M \models$ "$\mathbb{S}$ is $\mathbf{R}$-good" then, in $N$, $\mathbb{S}^N$ forces that every real in $X \cap N$ that is $\mathbf{R}$-unbounded over $M$ is still $\mathbf{R}$-unbounded over $M^{\mathbb{S}^M}$.*

**Corollary 3.12.** *Let $\Gamma \in M$ be a non-empty set. If $M \models$ "$\mathbb{B}_\Gamma$ is $\mathbf{R}$-good" then $\mathbb{B}_\Gamma^N$, in $N$, forces that every real in $X \cap N$ that is $\mathbf{R}$-unbounded over $M$ is still $\mathbf{R}$-unbounded over $M^{\mathbb{B}_\Gamma^M}$.*

**Lemma 3.13** ([BF11, Lemma 11], see also [Mej15, Lemma 5.13]). *Assume $\mathbb{P} \in M$ is a poset. Then, in $N$, $\mathbb{P}$ forces that every real in $X \cap N$ that is $\mathbf{R}$-unbounded over $M$ is still $\mathbf{R}$-unbounded over $M^{\mathbb{P}}$.*

**Lemma 3.14** (Blass and Shelah [BS89], [BF11, Lemmas 10, 12 and 13]). *Let $\mathbf{s}$ be a coherent pair of FS iterations (wlog $I^{\mathbf{s}} = \{0, 1\}$). Then, $\mathbb{P}_{0,\xi} \lessdot \mathbb{P}_{1,\xi}$ for all $\xi \le \pi$.*

*Moreover, if $\dot{c}$ is a $\mathbb{P}_{1,1}$-name of a real in $X$, $\pi$ is limit and $\mathbb{P}_{1,\xi}$ forces that $\dot{c}$ is $\mathbf{R}$-unbounded over $V_{0,\xi}$ for all $0 < \xi < \pi$, then $\mathbb{P}_{1,\pi}$ forces that $\dot{c}$ is $\mathbf{R}$-unbounded over $V_{0,\pi}$.*

We finish this part with the main result of this subsection.

**Theorem 3.15.** *Let $\mathbf{R} = \langle X, Y, \sqsubset \rangle$ be a Prs, $\theta$ a cardinal of uncountable cofinality and let $\mathbf{s}$ be a standard coherent system of FS iterations of length $\pi > 0$ that satisfies the hypothesis of Lemma 2.7. Further assume that*

(i) *$\Gamma \subseteq \Omega$ has size $\geq \theta$,*

(ii) *$D \in I$ and $\Gamma \subseteq D$,*

(iii) *for each $l \in \Gamma$, $\mathbb{P}_{D,1}$ adds a real $\dot{c}_l$ in $X$ such that, for any $A \in D$ in $I \cap [\Omega]^{<\theta}$, if $l \in D \smallsetminus A$ then $\mathbb{P}_{D,1}$ forces that $\dot{c}_l$ is $\mathbf{R}$-unbounded over $V_{A,1}$, and*

(iv) *for every $\xi \in S^{\mathbf{s}}$ and $B \in I \cap [\Omega]^{<\theta}$, $\mathbb{P}_{B,\xi}$ forces that $\dot{\mathbb{Q}}_{B,\xi}$ is $\mathbf{R}$-good.*

*Then $\mathbb{P}_{D,\pi}$ forces that the family $\dot{F} := \{\dot{c}_l : l \in \Gamma\}$ is strongly $\theta$-$\mathbf{R}$-unbounded. In particular, $\mathbb{P}_{D,\pi}$ forces $\mathfrak{b}(\mathbf{R}) \leq |\dot{F}|$ and, when $\theta$ is regular, this poset forces $|\dot{F}| \leq \mathfrak{d}(\mathbf{R})$.*

*Proof.* Let $\dot{y}$ be a $\mathbb{P}_{D,\pi}$-name of a member of $Y$. By Lemma 2.7, there is some $A \in I \cap [D]^{<\theta}$ such that $\dot{y}$ is a $\mathbb{P}_{A,\pi}$-name. Fix $l \in \Gamma \smallsetminus A$. By Lemmas 3.11, 3.13, 3.14 and Corollary 3.12 applied to the coherent pair $\mathbf{s} \restriction \{A, D\}$, $\mathbb{P}_{D,\pi}$ forces that $\dot{c}_l$ is $\mathbf{R}$-unbounded over $V_{A,\pi}$, which implies that $\dot{c}_l \not\sqsubset \dot{y}$. Therefore, $\mathbb{P}_{D,\pi}$ forces that $\{x \in \dot{F} : x \sqsubset \dot{y}\} \subseteq \{\dot{c}_l : l \in \Gamma \cap A\}$, which has size $\leq |A| < \theta$.

The second statement is a consequence of Remark 3.2                                  $\square$

### 3.2. Preservation of mad families.

**Definition 3.16.** Fix $A \subseteq [\omega]^{\aleph_0}$.

(1) Let $P \subseteq \left[ [\omega]^{\aleph_0} \right]^{<\aleph_0}$. For $x \subseteq \omega$ and $h : \omega \times P \to \omega$, define $x \sqsubset^* h$ by

$$\forall^\infty n < \omega \forall F \in P([n, h(n, F)) \smallsetminus \bigcup F \not\subseteq x).$$

(2) Define the relational system $\mathbf{Md}(A) := \langle [\omega]^{\aleph_0}, \omega^{\omega \times [A]^{<\aleph_0}}, \sqsubset^* \rangle$.

(3) Say that a poset $\mathbb{P}$ is *uniformly $\mathbf{Md}(A)$-good* if, for any $\mathbb{P}$-name $\dot{h}$ of a member of $\omega^{\omega \times [A]^{<\aleph_0}}$, there is a non-empty countable $H \subseteq \omega^{\omega \times [A]^{<\aleph_0}}$ (in the ground model) such that, for any countable $C \subseteq A$ and any $x \in [\omega]^{\aleph_0}$, if $x \not\sqsubset^* h' \restriction (\omega \times [C]^{<\aleph_0})$ for all $h' \in H$ then $\Vdash x \not\sqsubset^* \dot{h} \restriction (\omega \times [C]^{<\aleph_0})$.

Throughout this subsection, fix transitive models $M \subseteq N$ of ZFC and $A \in M$ such that $A \subseteq [\omega]^{\aleph_0} \cap M$. The relational system $\mathbf{Md}(A)$ helps to abbreviate the main notion presented in [BF11] for the preservation of mad families. What is defined in [BF11, Def. 2] as $(\star^{M,N}_{A,a})$ for $a \in [\omega]^{\aleph_0}$, which is the same as "$a$ diagonalizes $M$ outside $A$" in [FFMM18, Def. 4.2], actually means in our notation that $a$ is $\mathbf{Md}(A)$-unbounded over $M$. Note that, for any countable $C \subseteq [\omega]^{\aleph_0}$, $\mathbf{Md}(C)$ is a Prs.

The following results from [BF11] indicate that the a.d. family added by $\mathbb{H}_\Omega$ is composed of unbounded reals in the sense of relational systems like in Definition 3.16(2), which in turn becomes a mad family when $\Omega$ is uncountable.

**Lemma 3.17** ([BF11, Lemma 3]). *If $a^* \in [\omega]^{\aleph_0}$ is $\mathbf{Md}(A)$-unbounded over $M$ then $|a^* \cap x| = \aleph_0$ for any $x \in M \smallsetminus \mathcal{I}(A)$ where $\mathcal{I}(A) := \{x \subseteq \omega : \exists F \in [A]^{<\aleph_0}(x \subseteq^* \bigcup F)\}$.*

**Lemma 3.18** ([BF11, Lemma 4]). *Let $\Omega$ be a set, $z^* \in \Omega$ and $\dot{A} := \langle \dot{a}_z : z \in \Omega \rangle$ the a.d. family added by $\mathbb{H}_\Omega$. Then, $\mathbb{H}_\Omega$ forces that $\dot{a}_{z^*}$ is $\mathbf{Md}(\dot{A} \restriction (\Omega \smallsetminus \{z^*\}))$-unbounded over $V^{\mathbb{H}_{\Omega \smallsetminus \{z^*\}}}$.*

The known results about the preservation of $\mathbf{Md}(A)$-unbounded reals along coherent pairs of FS iterations are referred below. This is similar to the previous discussion about preservation of $\mathbf{R}$-unbounded reals for a Prs $\mathbf{R}$.

**Lemma 3.19** ([BF11, Lemma 12]). *Let $\mathbf{s}$ be a coherent pair of FS iterations (wlog $I^{\mathbf{s}} = \{0,1\}$) with $\pi = \pi^{\mathbf{s}}$ limit, $\dot{A}$ a $\mathbb{P}_{0,1}$-name of a family of infinite subsets of $\omega$ and $\dot{a}^*$ a $\mathbb{P}_{1,1}$-name for an infinite subset of $\omega$ such that*

$$\Vdash_{\mathbb{P}_{1,\xi}} \text{``}\dot{a}^* \text{ is } \mathbf{Md}(A)\text{-unbounded over } V_{0,\xi} \text{''}$$

*for all $0 < \xi < \pi$. Then, $\mathbb{P}_{0,\pi} \lessdot \mathbb{P}_{1,\pi}$ and $\Vdash_{\mathbb{P}_{1,\pi}} \text{``}\dot{a}^* \text{ is } \mathbf{Md}(A)\text{-unbounded over } V_{0,\pi} \text{''}$.*

**Lemma 3.20** ([BF11, Lemma 11]). *Let $\mathbb{P} \in M$ be a poset. If $N \models \text{``}a^* \text{ is } \mathbf{Md}(A)\text{-unbounded over } M\text{''}$ then*

$$N^{\mathbb{P}} \models \text{``}a^* \text{ is } \mathbf{Md}(A)\text{-unbounded over } M^{\mathbb{P}}\text{''}.$$

**Corollary 3.21.** *If $\Omega \in M$ and $N \models \text{``}a^* \text{ is } \mathbf{Md}(A)\text{-unbounded over } M\text{''}$ then*

$$N^{\mathbb{C}_\Omega} \models \text{``}a^* \text{ is } \mathbf{Md}(A)\text{-unbounded over } M^{\mathbb{C}_\Omega}\text{''}.$$

*Likewise, $\mathbb{H}_\Omega$ satisfies a similar statement.*

**Lemma 3.22** ([FFMM18, Lemma 4.8 and Cor. 4.11]). *Let $\mathbb{S}$ be either $\mathbb{E}$ or a random algebra. If $N \models \text{``}a^* \text{ is } \mathbf{Md}(A)\text{-unbounded over } M\text{''}$ then*

$$N^{\mathbb{S}^N} \models \text{``}a^* \text{ is } \mathbf{Md}(A)\text{-unbounded over } M^{\mathbb{S}^M}\text{''}.$$

The previous result indicates that $\mathbb{E}$ and random forcing, when used as iterands in a coherent pair of FS iterations, help to preserve $\mathbf{Md}(A)$-unbounded reals. To generalize this fact, we use the notion of 'uniformly good' introduced in Definition 3.16(3).

**Theorem 3.23.** *Let $\mathbb{S}$ be a Suslin ccc poset coded in $M$ and $A \in M$, $A \subseteq [\omega]^{\aleph_0}$. Assume*

  *($\star$) $[A]^{\aleph_0} \cap M$ is cofinal in $[A]^{\aleph_0} \cap N$.*

*If $M \models \text{``}\mathbb{S} \text{ is uniformly } \mathbf{Md}(A)\text{-good''}$ then, in $N$, $\mathbb{S}^N$ forces that every real in $[\omega]^{\aleph_0} \cap N$ that is $\mathbf{Md}(A)$-unbounded over $M$ is still $\mathbf{Md}(A)$-unbounded over $M^{\mathbb{S}^M}$.*

Note that, when $N$ is a generic extension of $M$ by a proper poset, ($\star$) holds.

*Proof.* Let $a \in [\omega]^{\aleph_0} \cap N$ be $\mathbf{Md}(A)$-unbounded over $M$. Assume that $\dot{h} \in M$ is a $\mathbb{S}^M$-name of a function in $\omega^{\omega \times [A]^{<\omega}}$. As $\mathbb{S}$ is uniformly $\mathbf{Md}(A)$-good in $M$, there is a family $\{h_n : n < \omega\} \subseteq \omega^{\omega \times [A]^{<\aleph_0}}$ (in $M$) that witnesses goodness for $\dot{h}$. Thus $a \not\sqsupset^* h_n$ for every $n < \omega$, so we can find a countable $C \subseteq A$ such that $a \not\sqsupset^* h_n{\upharpoonright}(\omega \times [C]^{<\aleph_0})$ for every $n < \omega$. By ($\star$), wlog we can find such $C$ in $M$.

In $M$ the statement "for every $x \in [\omega]^{\aleph_0}$, if $x \not\sqsupset^* h_n{\upharpoonright}(\omega \times [C]^{<\aleph_0})$ for all $n < \omega$, then $\Vdash_{\mathbb{S}} x \not\sqsupset^* \dot{h}{\upharpoonright}(\omega \times [C]^{<\aleph_0})$" is true. Furthermore, as this statement is a conjunction of a $\mathbf{\Sigma}^1_1$-statement with a $\mathbf{\Pi}^1_1$-statement of the reals (see e.g. [CM, Claim 4.27]), it is also true in $N$. In particular, since $a \not\sqsupset^* h_n{\upharpoonright}(\omega \times [C]^{<\aleph_0})$ for every $n < \omega$, $\Vdash^N_{\mathbb{S}^N} a \not\sqsupset^* \dot{h}{\upharpoonright}(\omega \times [C]^{<\aleph_0})$. $\square$

Though $\mathbb{E}$ and $\mathbb{B}$ are indeed uniformly $\mathbf{Md}(A)$-good (by Theorem 3.27 and Lemma 3.29), the application of Theorem 3.23 yields a version of Lemma 3.22 restricted to the condition ($\star$). To avoid this restriction, we consider an alternative generalization based on the following notion.

**Definition 3.24.** Let $\mathbb{P}$ be a poset.

(1) Say that a set $Q \subseteq \mathbb{P}$ is *Frechet-linked (in $\mathbb{P}$)*, abbreviated Fr-*linked*, if, for any sequence $\bar{p} = \langle p_n : n < \omega \rangle$ in $Q$, there is some $q \in \mathbb{P}$ that forces $\exists^\infty n < \omega(p_n \in \dot{G})$.

(2) Let $\mu$ be an infinite cardinal. Say that a poset $\mathbb{P}$ is *$\mu$-Frechet-linked* (often abbreviated $\mu$-Fr-*linked*) if there is a sequence $\langle Q_\alpha : \alpha < \mu \rangle$ of Fr-linked subsets of $\mathbb{P}$ such that $\bigcup_{\alpha < \mu} Q_\alpha$ is dense in $\mathbb{P}$.

By *$\sigma$-Fr-linked* we mean $\aleph_0$-Fr-linked.

(3) A poset $\mathbb{S}$ is *Suslin $\sigma$-Frechet-linked* if $\mathbb{S}$ is a subset of some Polish space, the relations $\leq$ and $\perp$ are $\mathbf{\Sigma}_1^1$ (in that Polish space) and $\mathbb{S} = \bigcup_{n < \omega} Q_n$ where each $Q_n$ is a Fr-linked $\mathbf{\Sigma}_1^1$ set.

Here, Fr denotes the Frechet filter on $\omega$. The reason of the terminology 'Frechet-linked' is that this notion corresponds to a particular case on Fr of a more general notion of linkedness with filters that we provide in Example 5.4.

**Remark 3.25.** (1) The notion '$\mu$-Fr-linked' is a forcing property, i.e., if $\mathbb{P}$ and $\mathbb{Q}$ are posets, $\mathbb{P} \lessdot \mathbb{Q}$ (in the sense that the Boolean completion of $\mathbb{P}$ is completely embedded into the completion of $\mathbb{Q}$) and $\mathbb{Q}$ is $\mu$-Fr-linked then so is $\mathbb{P}$ (see more on this in Section 5).

(2) No Fr-linked subset of a poset can contain infinite antichains. In addition, if $\mathbb{P}$ is a poset and $Q \subseteq \mathbb{P}$, the statement "$Q$ does not contain infinite antichains" is absolute for transitive models of ZFC. This is because that statement is equivalent to say that "$T$ is a well-founded tree" where $T := \{s \in Q^{<\omega} : \mathrm{ran}\, s$ is an antichain$\}$.

(3) As a consequence of (2), $\mu$-Fr-linked posets are $\mu^+$-cc. Even more, by [HT48, Thm. 2.4], they are $\mu^+$-Knaster (see more in Section 5).

(4) Any poset of size $\leq \mu$ is $\mu$-Fr-linked (witnessed by its singletons). In particular, Cohen forcing is $\sigma$-Fr-linked.

(5) By (3), any Suslin $\sigma$-Fr-linked poset is Suslin ccc. Moreover, if $\langle Q_n : n < \omega \rangle$ witnesses that a poset $\mathbb{S}$ is Suslin $\sigma$-Fr-linked then the statement "$Q_n$ is Fr-linked" is $\mathbf{\Pi}_2^1$ (by (6) below, its negation is equivalent to $\exists f \in Q_n^\omega \exists g \in \mathbb{S}^\omega(\{g(n) : n < \omega\}$ is a maximal antichain and $\forall n < \omega \exists m < \omega \forall k \geq m(g(n) \perp f(k)))$, which is $\mathbf{\Sigma}_2^1$). Therefore, if $M \models$"$\mathbb{S}$ is Suslin $\sigma$-Fr-linked" and $\omega_1^N \subseteq M$ then $N \models$"$\mathbb{S}$ is Suslin $\sigma$-Fr-linked".

(6) Let $\mathbb{P}$ be a poset and $Q \subseteq \mathbb{P}$. Note that a sequence $\langle p_n : n < \omega \rangle$ in $Q$ witnesses that $Q$ is <u>not</u> Fr-linked iff the set

$$\{q \in \mathbb{P} : \forall^\infty n < \omega(q \perp p_n)\}$$

is dense.

**Lemma 3.26.** *Let $\mathbb{P}$ be a poset and $Q \subseteq \mathbb{P}$ Fr-linked. If $\dot{n}$ is a $\mathbb{P}$-name of a natural number then there is some $m < \omega$ such that $\forall p \in Q(p \nVdash m < \dot{n})$.*

*Proof.* Towards a contradiction, assume that for each $m < \omega$ there is some $p_m \in Q$ that forces $m < \dot{n}$. Hence, as $Q$ is Fr-linked, there is some $q \in Q$ that forces $\exists^\infty m < \omega(p_m \in \dot{G})$, which implies $\exists^\infty m < \omega(m < \dot{n})$, a contradiction. $\square$

**Theorem 3.27.** *Any $\sigma$-Fr-linked poset is uniformly $\mathbf{Md}(B)$-good for any $B \subseteq [\omega]^{\aleph_0}$.*

*Proof.* Let $\mathbb{P}$ be a $\sigma$-Fr-linked poset witnessed by $\langle Q_n : n < \omega \rangle$. Assume that $\dot{h}$ is a $\mathbb{P}$-name of a function in $\omega^{\omega \times [B]^{<\omega}}$. Fix $n < \omega$. For each $k < \omega$ and $F \in [B]^{<\aleph_0}$, by Lemma 3.26 there is some $h_n(k, F) < \omega$ such that $\forall p \in Q_n(p \nVdash h_n(k, F) < \dot{h}(k, F))$. This allows to define $h_n \in \omega^{\omega \times [B]^{<\omega}}$.

Assume that $C \subseteq B$ is infinite and that $x \in [\omega]^{\aleph_0}$ is $\mathbf{Md}(C)$-unbounded over $\{h_n \upharpoonright (\omega \times [C]^{<\aleph_0}) : n < \omega\}$. We show that $\Vdash x \not\sqsubseteq^* \dot{h} \upharpoonright (\omega \times [C]^{<\aleph_0})$. Assume that $p \in \mathbb{P}$ and $m < \omega$. Wlog, we may assume that $p \in Q_n$ for some $n < \omega$. As $x \not\sqsubseteq^* h_n \upharpoonright (\omega \times [C]^{<\aleph_0})$, there are some $k > m$ and $F \in [C]^{<\aleph_0}$ such that $[k, h_n(k, F)) \smallsetminus \bigcup F \subseteq x$. On the other hand, by the definition of $h_n$, there is some $q \leq p$ in $\mathbb{P}$ that forces $\dot{h}(k, F) \leq h_n(k, F)$. Hence $q \Vdash [k, \dot{h}(k, F)) \smallsetminus \bigcup F \subseteq x$. □

As a consequence, Theorem 3.23 can be applied to Suslin $\sigma$-Fr-linked posets. However, it can be proved directly that Suslin $\sigma$-Fr-linked posets preserves $\mathbf{Md}(A)$-unbounded reals without the condition (⋆).

**Theorem 3.28.** *Let $\mathbb{S}$ be a Suslin $\sigma$-Fr-linked poset coded in $M$. Then, in $N$, $\mathbb{S}^N$ forces that every real in $[\omega]^{\aleph_0} \cap N$ that is $\mathbf{Md}(A)$-unbounded over $M$ is still $\mathbf{Md}(A)$-unbounded over $M^{\mathbb{S}^M}$.*

*Proof.* Let $\langle Q_n : n < \omega \rangle \in M$ be a sequence that witnesses that $\mathbb{S}$ is Suslin $\sigma$-Fr-linked. Let $a \in [\omega]^{\aleph_0} \cap N$ be $\mathbf{Md}(A)$-unbounded over $M$. Assume that $\dot{h} \in M$ is a $\mathbb{S}^M$-name of a function in $\omega^{\omega \times [A]^{<\omega}}$. Fix $p \in \mathbb{S}^N$ and $m < \omega$, so there is some $n < \omega$ such that $p \in Q_n^N$. Down in $M$, find $h_n \in \omega^{\omega \times [A]^{<\omega}} \cap M$ as in the proof of Theorem 3.27. As $a \not\sqsubseteq^* h_n$, there are some $k \geq m$ and $F \in [A]^{<\aleph_0}$ such that $[k, h_n(k, F)) \smallsetminus \bigcup F \subseteq a$. On the other hand, the statement "no $q \in Q_n$ forces that $h_n(k, F) < \dot{h}(k, F)$" is $\mathbf{\Pi}_1^1$, so it is absolute and, as true in $M$, it also holds in $N$. Therefore, in $N$, $p$ does not force $h_n(k, F) < \dot{h}(k, F)$, which implies that some $q \leq p$ in $\mathbb{S}^N$ forces the contrary. This clearly implies that $q$ forces $[k, \dot{h}(k, F)) \smallsetminus \bigcup F \subseteq a$. □

Therefore, in conjunction with the following result, the previous theorem is a suitable generalization of Lemma 3.22.

**Lemma 3.29.** *The posets $\mathbb{E}$ and $\mathbb{B}$ are Suslin $\sigma$-Fr-linked. Moreover, any complete Boolean algebra that admits a strictly positive $\sigma$-additive measure (e.g. any random algebra) is $\sigma$-Fr-linked.*

*Proof.* For each $s \in \omega^{<\omega}$ and $m < \omega$ define $E_{s,m} := \{(t, \varphi) \in \mathbb{E} : t = s$ and $\forall i < \omega (|\varphi(i)| \leq m)\}$. This set is actually Borel in $\omega^{<\omega} \times \mathcal{P}(\omega)^\omega$ (the Polish space where $\mathbb{E}$ is defined). A compactness argument similar to the one in [Mil81] shows that $E_{s,n}$ is Fr-linked in $\mathbb{E}$. Fix a non-principal ultrafilter $U$ on $\omega$ and let $\langle p_n : n < \omega \rangle$ be a sequence in $E_{s,m}$. Write $p_n = (s, \varphi_n)$. For each $i < \omega$ define $\varphi(i) \subseteq \omega$ such that $l \in \varphi(i)$ iff $\{n < \omega : l \in \varphi_n(i)\} \in U$. It can be proved that $|\varphi(i)| \leq m$, (so $q := (s, \varphi) \in \mathbb{E}$) and that $q$ forces $\exists^\infty n < \omega (p_n \in \dot{G})$ (that is, for any $q' \leq q$ and $n < \omega$, there is some $k \geq n$ such that $q'$ is compatible with $p_k$).

Now, consider random forcing as $\mathbb{B} = \bigcup_{m < \omega} B_m$ where[1]

$$B_m := \left\{ T \subseteq 2^{<\omega} : T \text{ is a well-pruned tree and } \mathrm{Lb}([T]) \geq \frac{1}{m+1} \right\}$$

Note that $B_m$ is Borel in $2^{2^{<\omega}}$. It is enough to show that $B_m$ is Fr-linked. Assume the contrary, so by Remark 3.25(6) there are a sequence $\langle T_n : n < \omega \rangle$ in $B_m$ and a partition $\langle A_n : n < \omega \rangle$ of $2^\omega$ into Borel sets of positive measure such that, for each $n < \omega$, $A_n \cap [T_k]$ has measure zero for all but finitely many $k < \omega$. Construct an increasing function $g : \omega \to \omega$ such that $A_n \cap [T_k]$ has measure zero for all $k \geq g(n)$. As $2^\omega = \bigcup_{n < \omega} A_n$, we can find some $n^* < \omega$ such that the measure of $A^* := \bigcup_{n < n^*} A_n$ is strictly larger than $1 - \frac{1}{m+1}$. Hence $A^* \cap [T]$ has positive measure for any $T \in B_m$, but this contradicts that $A^* \cap [T_k]$ has measure zero for all $k \geq g(n^* - 1)$.

---

[1] A *well-pruned tree* is a non-empty tree such that every node has a successor.

A similar proof works for any complete Boolean algebra that admits a strictly positive $\sigma$-additive measure.                                                                         $\square$

Indeed, the notion of $\mu$-Fr-linked behaves well for preservation of $\mathbf{D}$-unbounded families as shown in the following result. Even more, this generalizes the facts of Example 3.7(2).

**Theorem 3.30.** *If $\mu < \theta$ are infinite cardinals then any $\mu$-Fr-linked poset is $\theta$-$\mathbf{D}$-good.*

*Proof.* This argument is very similar to the proof of Theorem 3.27. Let $\mathbb{P}$ be a $\mu$-Fr-linked poset witnessed by $\langle Q_\alpha : \alpha < \mu \rangle$, and let $\dot{y}$ be a $\mathbb{P}$-name of a member of $\omega^\omega$. Using Lemma 3.26, for each $\alpha < \mu$ find a $y_\alpha \in \omega^\omega$ such that, for every $i < \omega$, no member of $Q_\alpha$ forces that $y_\alpha(i) < \dot{y}(i)$. It can be proved that $\{y_\alpha : \alpha < \mu\}$ witnesses goodness for $\dot{y}$.   $\square$

**Remark 3.31.** There is a Suslin $\sigma$-Fr-linked poset that is not $\mathbf{Lc}$-good. For $b, h \in \omega^\omega$ such that $\forall i < \omega(b(i) > 0)$ and $h$ goes to infinity, consider the poset

$$\mathbb{LOC}_{b,h} := \left\{ p \in \prod_{i < \omega} [b(i)]^{\leq h(i)} : \exists m < \omega \forall^\infty i < \omega(|p(i)| \leq m) \right\}$$

ordered by $q \leq p$ iff $\forall i < \omega(p(i) \subseteq q(i))$. This poset adds a slalom $\dot{\varphi}_* \in \prod_{i < \omega} [b(i)]^{\leq h(i)}$, defined by $\dot{\varphi}_*(i) := \bigcup_{p \in \dot{G}} p(i)$, such that $x \in^* \varphi_*$ for every $x \in \prod_{i < \omega} b(i)$ in the ground model. Note that $\dot{\varphi}_*(i)$ is forced to either have size $h(i)$ (whenever $h(i) \leq b(i)$) or to be equal to $b(i)$ (whenever $b(i) \leq h(i)$).

This poset is Suslin $\sigma$-Fr-linked. In fact, for any $s \in \bigcup_{n < \omega} \prod_{i < n} b(i)$ and $m < \omega$, the set

$$L_{b,h}(s, m) := \left\{ p \in \prod_{i < \omega} [b(i)]^{\leq h(i)} : s \subseteq p, \text{ and } \forall i \geq |s|(|p(i)| \leq m) \right\}$$

is Borel and Fr-linked, and $\mathbb{LOC}_{b,h} = \bigcup_{s,m} L_{b,h}(s, m)$. To see that $L_{b,h}(s, m)$ is Fr-linked, assume that $\langle p_n : n < \omega \rangle$ is a sequence in $L_{b,h}(s, m)$ and fix a non-principal ultrafilter $U$ on $\omega$. For each $i \geq |s|$, as $[b(i)]^{\leq m}$ is finite, we can find $q(i) \in [b(i)]^{\leq m}$ and $a_i \in U$ such that $p_n(i) = q(i)$ for any $n \in a_i$. Hence $q \in L_{b,h}(s, m)$ where $q(i) := s(i)$ for all $i < |s|$. It remains to show that $q$ forces $|\{n < \omega : p_n \in \dot{G}\}| = \aleph_0$. If $r \leq q$ and $n_0 < \omega$, then we can find $k, m_0 < \omega$ such that $k \geq |s|$ and, for any $i \geq k$, $|r(i)| \leq m_0$ and $m_0 + m \leq h(i)$. Choose some $n \in \bigcap_{i < k} a_i \smallsetminus n_0$ (put $a_i := \omega$ for $i < |s|$). Hence $q{\restriction}k = p_n{\restriction}k$ and $q'$ forces that $p_n \in \dot{G}$ where $q'(i) := r(i) \cup p_n(i)$ (for $i < k$, $q'(i) = r(i)$; for $i \geq k$, $|q'(i)| \leq m_0 + m \leq h(i)$).

Now, if $h = \mathrm{id}$ and $b(i) = i^i$ for every $i < \omega$ then $\mathbb{LOC}_{b,h}$ is <u>not</u> $\mathbf{Lc}$-good. Consider the name $\dot{\varphi}_*$ of the slalom that this poset adds, which is clearly a name of a member of $\mathcal{S}(\omega, \mathcal{H})$ (see Example 3.7(4)). If $H$ is a countable subset of $\mathcal{S}(\omega, \mathcal{H})$, then there exists an $x \in \prod_{i < \omega} b(i)$ such that $x$ is not localized by any member of $H$. On the other hand, $\mathbb{LOC}_{b,h}$ already forces that $x \in^* \dot{\varphi}_*$.

We finish this section with a general result about preservation of mad families through standard coherent systems.

**Theorem 3.32.** *Let $\nu$ be a cardinal of uncountable cofinality and let $\mathbf{s}$ be a standard coherent system that satisfies the hypothesis of Lemma 2.7 for $\theta = \nu$. Further assume that*

(i) $\Gamma \subseteq \Omega$ has size $\geq \nu$,

(ii) $D \in I$ and $\Gamma \subseteq D$,

(iii) for each $l \in \Gamma$, $\mathbb{P}_{D,1}$ adds a real $\dot{a}_l$ in $[\omega]^{\aleph_0}$ such that, for any $Z \subseteq D$ in $I \cap [\Omega]^{<\nu}$, $\dot{A}{\restriction}Z := \langle \dot{a}_l : l \in Z \cap \Gamma \rangle$ is a $\mathbb{P}_{Z,1}$-name and, whenever $l \in D \smallsetminus Z$, $\mathbb{P}_{D,1}$ forces that $\dot{a}_l$ is $\mathbf{Md}(\dot{A}{\restriction}Z)$-unbounded over $V_{Z,1}$, and

(iv) for every $\xi \in S^{\mathbf{s}}$ and $B \in I \cap [D]^{<\nu}$, $\mathbb{P}_{B,\xi}$ forces that $\dot{\mathbb{Q}}_{B,\xi}$ is either uniformly $\mathbf{Md}(\dot{A}{\restriction}B)$-good or a random algebra.

Then $\mathbb{P}_{D,\pi}$ forces that any infinite subset of $\omega$ intersects some member of $\dot{A} := \dot{A}{\restriction}\Gamma$. In particular, $\mathbb{P}_{D,\pi}$ forces $\mathfrak{a} \leq |\Gamma|$ whenever $\dot{A}$ is an a.d. family.

*Proof.* Let $\dot{x}$ be a $\mathbb{P}_{D,\pi}$-name of an infinite subset of $\omega$. By Lemma 2.7, there is some $Z \subseteq D$ in $I \cap [\Omega]^{<\nu}$ such that $\dot{x}$ is a $\mathbb{P}_{Z,\pi}$-name. Thus, by Lemmas 3.19, 3.20 and Theorem 3.23, for any $l \in \Gamma \smallsetminus Z$, $\mathbb{P}_{D,\pi}$ forces that $\dot{a}_l$ is $\mathbf{Md}(\dot{A}{\restriction}Z)$-unbounded over $V_{Z,\pi}$, so $\dot{x} \notin \mathcal{I}(\dot{A}{\restriction}Z)$ implies that $\dot{x} \cap \dot{a}_l$ is infinite by Lemma 3.17. As $\dot{x} \in \mathcal{I}(\dot{A}{\restriction}Z)$ implies that $\dot{x} \cap \dot{a}_j$ is infinite for some $j \in Z \cap \Gamma$, we are done. $\qquad\square$

**Remark 3.33.** By Theorem 3.27, it is clear that in condition (iv) of Theorem 3.32 we can use Suslin $\sigma$-Fr-linked posets.

## 4. Applications

The following result generalizes [FFMM18, Thm. 4.17] in the sense that it allows to preserve mad families of singular cardinality along a more general type of FS iterations.

**Theorem 4.1.** *Let $\theta$ be an uncountable regular cardinal and let $\Omega$ be a set of size $\geq \theta$. After forcing with $\mathbb{H}_\Omega$, any further FS iteration where each iterand is one of the following types preserves the mad family added by $\mathbb{H}_\Omega$.*

*(0) Suslin $\sigma$-Fr-linked.*

*(1) Random algebra.*

*(2) Hechler poset (for adding a mad family).*

*(3) Poset with ccc of size $< \theta$.*

*Proof.* Consider a FS iteration $\mathbb{P}_\pi = \langle \mathbb{P}_\xi, \dot{\mathbb{Q}}_\xi : \xi < \pi \rangle$ with $\pi > 0$ such that $\dot{\mathbb{Q}}_0 = \mathbb{P}_1 = \mathbb{H}_\Omega$ and, for $0 < \xi < \pi$, $\dot{\mathbb{Q}}_\xi$ is of one of the types above. To be more precise, let $\langle C_j : j < 4 \rangle$ be a partition of $[1, \pi)$ such that, for each $j < 4$ and $\xi \in C_j$, $\dot{\mathbb{Q}}_\xi$ is a $\mathbb{P}_\alpha$-name of a poset of type $(j)$. Note that this iteration can be defined as the standard coherent system $\mathbf{m}$ on $I^{\mathbf{m}} := \langle \mathcal{P}(\Omega), \subseteq \rangle$ such that

(o) $\mathbb{P}_{X,1} = \mathbb{H}_X$ for any $X \subseteq \Omega$;

(i) $S^{\mathbf{m}} = C_0 \cup C_1$, $C^{\mathbf{m}} = C_2 \cup C_3$;

(ii) $\Delta : [1, \pi) \to [\Omega]^{<\theta}$ such that $\Delta(\xi) = \emptyset$ for each $\xi \in C_1 \cup C_2$;

(iii) for $\xi \in C_0$, when $\dot{\mathbb{Q}}_\xi$ is coded in $V_{\Delta(\xi),\xi}$, $\dot{\mathbb{S}}_\xi = \dot{\mathbb{Q}}_\xi$ (or trivial otherwise, though this latter case will not happen) and, for $\xi \in C_1$, $\mathbb{S}_\xi$ is the random algebra $\dot{\mathbb{Q}}_\xi$ (wlog its support is in the ground model);

(iv) for $\xi \in C_2$, $\dot{\mathbb{Q}}_\xi^{\mathbf{m}} = \dot{\mathbb{Q}}_\xi$ (wlog, the support of this Hechler poset is in the ground model) and, for $\xi \in C_3$, when $\dot{\mathbb{Q}}_\xi$ is forced to be in $V_{\Delta(\xi),\xi}$, $\dot{\mathbb{Q}}_\xi^{\mathbf{m}} = \dot{\mathbb{Q}}_\xi$ (or trivial otherwise, though this will not happen).

By recursion on $\xi \leq \pi$, $\mathbf{m}{\restriction}\xi$ and $\Delta{\restriction}\xi$ should be constructed and it should be guaranteed that $\mathbb{P}_{\Omega,\xi} = \mathbb{P}_\xi$. This is fine in the steps $\xi = 0, 1$ and the limit steps. Consider the successor step, i.e., assume we have the construction up to $\xi$. As $\xi \in C_j$ for some $j < 4$, we consider cases for each $j$. If $j = 0$ then, as a Suslin ccc poset is coded by reals and $\mathbb{P}_{\Omega,\xi} = \mathbb{P}_\xi$, by Lemma 2.7 there is some $\Delta(\xi) \in [\Omega]^{<\theta}$ such that $\dot{Q}_\xi$ is (coded by) a $\mathbb{P}_{\Delta(\xi),\xi}$-name; if $j = 1, 2$ then, as the support of $\dot{Q}_\xi$ can be assumed to be in the ground model, we can put $\Delta(\xi) = 0$; if $j = 3$ then, by Lemma 2.7 and the regularity of $\theta$, there is some $\Delta(\xi) \in [\Omega]^{<\theta}$ such that $\dot{Q}_\xi$ is a $\mathbb{P}_{\Delta(\xi),\xi}$-name. Therefore, in any case, we can define $\mathbf{m}{\upharpoonright}(\xi + 1)$ as required and it is clear that $\mathbb{P}_{\Omega,\xi} = \mathbb{P}_\xi$.

Now let $\dot{A} = \langle \dot{a}_l : l \in \Omega \rangle$ be the $\mathbb{P}_{\Omega,1}$-name of the generic a.d. family it adds. Note that $\dot{a}_l$ is a $\mathbb{P}_{\{l\},1}$-name and, by Lemma 3.18, for any $B \subseteq \Omega$ with $l \notin B$, $\mathbb{P}_{B \cup \{l\},1}$ forces that $\dot{a}_l$ is $\mathbf{Md}(\dot{A}{\upharpoonright}B)$-unbounded over $V_{B,1}$. Hence, by Theorem 3.32, $\mathbb{P}_{\Gamma,\pi}$ forces that $\dot{A}$ is a mad family.                                                                                              $\square$

**Remark 4.2.** The previous theorem remains true if we add the type

(0') Suslin ccc poset coded in the ground model such that, for any $X \subseteq \Omega$ in the ground model, it is uniformly $\mathbf{Md}(\dot{A}{\upharpoonright}X)$-good in any ccc generic extension of $V^{\mathbb{H}_X}$.

By Theorem 3.27 Suslin $\sigma$-Fr-linked posets coded in the ground model satisfy (0').

The remaining results in this section are improvements of the consistency results of [FFMM18, Sect. 5] about separating cardinals in Cichoń's diagram. Not only can we force an additional singular value, but the constructions are uniform in the sense that there is no need to distinguish between 2D or 3D constructions anymore since all the coherent systems can be constructed on a partial order of the form $\langle \mathcal{P}(\Omega), \subseteq \rangle$. In the following proofs, sum and product denote the corresponding operations in the ordinals, even when they are applied to cardinal numbers.

The following result improves [FFMM18, Thm. 5.2] about separating Cichoń's diagram into 7 different values.

**Theorem 4.3.** *Assume that $\theta_0 \leq \theta_1 \leq \kappa \leq \mu$ are uncountable regular cardinals, $\nu \leq \lambda$ are cardinals such that $\mu \leq \nu$, $\nu^{<\kappa} = \nu$ and $\lambda^{<\theta_1} = \lambda$. Then there is a ccc poset that forces* $\mathrm{MA}_{<\theta_0}$, $\mathrm{add}(\mathcal{N}) = \theta_0$, $\mathrm{cov}(\mathcal{N}) = \theta_1$, $\mathfrak{b} = \mathfrak{a} = \kappa$, $\mathrm{non}(\mathcal{M}) = \mathrm{cov}(\mathcal{M}) = \mu$, $\mathfrak{d} = \nu$ *and* $\mathrm{non}(\mathcal{N}) = \mathfrak{c} = \lambda$.

*Proof.* Let $\Omega_0$ and $\Omega_1$ be disjoint sets of size $\kappa$ and $\nu$ respectively. Put $\Omega := \Omega_0 \cup \Omega_1$. As $\nu^{<\kappa} = \nu$, we can enumerate $[\Omega]^{<\kappa} := \{ W_\zeta : \zeta < \nu \}$. Fix a bijection $g = (g_0, g_1, g_2) : \lambda \to 2 \times \nu \times \lambda$ and a function $t : \nu\mu \to \nu$ such that, for any $\zeta < \nu$, $t^{-1}[\{\zeta\}]$ is cofinal in $\nu\mu$ [2]. Put $\pi := \lambda\nu\mu$, $S := \{ \lambda\rho : \rho < \nu\mu \}$ and define $\Delta : [1, \pi) \to [\Omega]^{<\kappa}$ such that $\Delta(\lambda\rho) := \emptyset$, $\Delta(\lambda\rho + 1) := W_{t(\rho)}$ and $\Delta(\lambda\rho + 2 + \varepsilon) = W_{g_1(\varepsilon)}$ for each $\rho < \nu\mu$ and $\varepsilon < \lambda$.

Define the standard coherent system $\mathbf{m}$ of FS iterations of length $\pi$ on $\langle \mathcal{P}(\Omega), \subseteq \rangle$ such that $S^{\mathbf{m}} := S$, $C^{\mathbf{m}} := [1, \pi) \smallsetminus S$, $\Delta^{\mathbf{m}} := \Delta$, $\dot{\mathbb{Q}}^{\mathbf{m}}_{X,0} := \mathbb{H}_{X \cap \Omega_0} \times \mathbb{C}_{X \cap \Omega_1}$ and where the FS iterations at each interval of the form $[\lambda\rho, \lambda(\rho+1))$ for $\rho < \mu\nu$ is defined as follows. Assume that $\mathbf{m}{\upharpoonright}\lambda\rho$ has already been defined. For each $\zeta < \nu$ choose

(0) an enumeration $\{ \dot{Q}_{0,\zeta,\gamma} : \gamma < \lambda \}$ of all the (nice) $\mathbb{P}_{W_\zeta,\lambda\rho}$-names for posets of size $< \theta_0$, with underlying set contained in $\theta_0$, that are forced by $\mathbb{P}_{\Omega,\lambda\rho}$ to have ccc; and

(1) an enumeration $\{ \dot{Q}_{1,\zeta,\gamma} : \gamma < \lambda \}$ of all the (nice) $\mathbb{P}_{W_\zeta,\lambda\rho}$-names for subalgebras of random forcing of size $< \theta_1$.

Put $\mathbb{S}^{\mathbf{m}}_{\lambda\rho} := \mathbb{E}$ (only when $\rho > 0$) and, for each $\xi \in (\lambda\rho, \lambda(\rho + 1))$, put

(i) $\dot{\mathbb{Q}}^{\mathbf{m}}_\xi := \mathbb{D}^{V_{\Delta(\xi),\xi}}$ when $\xi = \lambda\rho + 1$, and

(ii) $\dot{\mathbb{Q}}^{\mathbf{m}}_\xi := \dot{\mathbb{Q}}_{g(\varepsilon)}$ when $\xi = \lambda\rho + 2 + \varepsilon$ for some $\varepsilon < \lambda$.

---

[2]For example, define $t(\nu\delta + \alpha) = \alpha$ for each $\delta < \mu$ and $\alpha < \nu$

This construction is possible because, as $\lambda^{<\theta_1} = \lambda$, each $\mathbb{P}_{\Omega,\xi}$ has size $\leq \lambda$.

It remains to show that $\mathbb{P} := \mathbb{P}_{\Omega,\pi}$ forces what we want. First note that this poset can be obtained by the FS iteration $\langle \mathbb{P}_{\Omega,\xi}, \dot{\mathbb{Q}}_{\Omega,\xi} : \xi < \pi \rangle$, and observe that all these iterands are $\theta_0$-**Lc**-good and $\theta_1$-**Cn**-good. Hence, by Theorem 3.10, $\mathbb{P}$ forces $\mathrm{add}(\mathcal{N}) \leq \theta_0$, $\mathrm{cov}(\mathcal{N}) \leq \theta_1$ and $\lambda \leq \mathrm{non}(\mathcal{N})$. Actually, those are equalities, even more, $\mathbb{P}$ forces $\mathrm{MA}_{<\theta_0}$ (which implies $\mathrm{add}(\mathcal{N}) \geq \theta_0$). To see this, let $\dot{\mathbb{R}}$ is a $\mathbb{P}$-name of a ccc poset of size $< \theta_0$ and $\dot{\mathcal{D}}$ a family of size $< \theta_0$ of dense subsets of $\dot{\mathbb{R}}$. By Lemma 2.7 there is some $\zeta < \nu$ such that both $\dot{\mathbb{R}}$ and $\dot{\mathcal{D}}$ are $\mathbb{P}_{W_\zeta,\pi}$-names. Moreover, as $\pi$ has cofinality $\mu$, there is some $\rho < \nu\mu$ such both are $\mathbb{P}_{W_\zeta,\lambda\rho}$-names. Therefore, there is some $\gamma < \lambda$ such that $\dot{\mathbb{R}} = \dot{\mathbb{Q}}_{0,\zeta,\gamma}$, so the generic set added by $\dot{\mathbb{Q}}_{g(\varepsilon)} = \dot{\mathbb{Q}}^{\mathbf{m}}_{W_\zeta,\xi}$ intersects all the dense sets in $\dot{\mathcal{D}}$ where $\varepsilon := g^{-1}(0,\zeta,\gamma)$ and $\xi := \lambda\rho + 2 + \varepsilon$. In a similar way, it can be proved that $\mathbb{P}$ forces $\mathrm{cov}(\mathcal{N}) \geq \theta_1$. On the other hand, since $\Vdash_\mathbb{P} \mathfrak{c} \leq \lambda$ follows from $|\mathbb{P}| \leq \lambda$, together with $\mathrm{non}(\mathcal{N}) \geq \lambda$ (see above) it is forced that $\mathrm{non}(\mathcal{N}) = \mathfrak{c} = \lambda$.

As the FS iteration that determines $\mathbb{P}$ has cofinality $\mu$ and $\mu$-cofinally many full eventually different reals are added by $\mathbb{E}$, $\mathbb{P}$ forces $\mathrm{cov}(\mathcal{M}) \leq \mu \leq \mathrm{non}(\mathcal{M})$. Actually, $\mathrm{non}(\mathcal{M}) \leq \mu \leq \mathrm{cov}(\mathcal{M})$ is forced by Theorem 3.10 applied to the Prs **Ed**, so $\mathrm{non}(\mathcal{M}) = \mathrm{cov}(\mathcal{M}) = \mu$. Now we show that $\mathbb{P}$ forces $\mathfrak{a} \leq \kappa$ and $\nu \leq \mathfrak{d}$. Let $\dot{A} := \{\dot{a}_l : l \in \Omega_0\}$ be the $\mathbb{H}_{\Omega_0}$-name of the mad family added by $\mathbb{H}_{\Omega_0}$ and let $\{\dot{c}_l : l \in \Omega_1\} \subseteq \omega^\omega$ be the Cohen reals added by $\mathbb{C}_{\Omega_1}$. For any $X \subseteq \Omega$, $l \in \Omega_0$ and $l' \in \Omega_1$, it is clear that $\dot{a}_l$ is a $\mathbb{P}_{X,1}$-name whenever $l \in X$, and $\dot{c}_{l'}$ is a $\mathbb{P}_{X,1}$-name whenever $l' \in X$. On the other hand, if $l' \notin X$ then $\mathbb{P}_{X \cup \{l'\},1}$ forces that $\dot{c}_{l'}$ is Cohen over $V_{X,1}$, hence it is **D**-unbounded over it; and if $l \notin X$ then $\mathbb{P}_{X \cup \{l\},1}$ forces that $\dot{a}_l$ is $\mathbf{Md}(\dot{A} \restriction (X \cap \Omega_0))$-unbounded over $V_{X,1}$. The latter is a consequence of Lemma 3.20 applied to $\mathbb{C}_{X \cap \Omega_1}$. Therefore, by Theorems 3.15 and 3.32 applied to $\{\dot{c}_{l'} : l' \in \Omega_1\}$ and $\{\dot{a}_l : l \in \Omega_0\}$ respectively, $\mathbb{P}$ forces $\nu \leq \mathfrak{d}$ (because $\{\dot{c}_{l'} : l' \in \Omega_1\}$ is strongly $\kappa$-**D**-unbounded) and $\mathfrak{a} \leq \kappa$.

It remains to show that $\mathbb{P}$ forces $\kappa \leq \mathfrak{b}$ and $\mathfrak{d} \leq \nu$. For each $\rho < \nu\mu$ denote by $\dot{d}_\rho$ the (restricted) dominating real over $V_{W_{t(\rho)},\lambda\rho+1}$ added by $\dot{\mathbb{Q}}_{W_{t(\rho)},\lambda\rho+1}$. It is enough to show that $\mathbb{P}$ forces that any subset of $\omega^\omega$ of size $< \kappa$ is dominated by some $\dot{d}_\rho$ (hence $\{\dot{d}_\rho : \rho < \nu\mu\}$ is a dominating family of size $\nu$). Let $\dot{F}$ be a $\mathbb{P}$-name of such a subset of $\omega^\omega$. By Lemma 2.7 and because $\mathrm{cf}(\pi) = \mu$, there are $\zeta < \nu$ and $\rho_0 < \nu\mu$ such that $\dot{F}$ is a $\mathbb{P}_{W_\zeta,\lambda\rho_0}$-name. Thus, there is some $\rho \in [\rho_0, \nu\mu)$ such that $t(\rho) = \zeta$, so $\mathbb{P}_{W_\zeta,\lambda\rho+2}$ forces that $\dot{d}_\rho$ dominates $\dot{F}$. $\qquad\square$

We summarize in the rest of this section the results from [Mej13a] and [FFMM18, Sect. 5] that can be improved by the method of the previous proof. Note that, in the forcing constructions for Theorems 4.5(b) and 4.6(c),(d), we cannot preserve a mad family added by a poset of the form $\mathbb{H}_\Omega$ because their constructions require that full generic dominating reals are added. For these items, it is enough to base their constructions on $\langle \mathcal{P}(\nu), \subseteq \rangle$ and start with $\dot{\mathbb{Q}}_{X,0} := \mathbb{C}_X$ for any $X \subseteq \nu$. In addition, by an argument similar to [FFMM18, Rem. 5.9], it can be additionally forced within these items that $\mathfrak{a} = \mu$.

**Theorem 4.4.** *Let* $\theta_0 \leq \theta_1 \leq \kappa \leq \nu \leq \lambda$ *be as in the statement of Theorem 4.3. Then there is a ccc poset forcing* $\mathrm{MA}_{<\theta_0}$, $\mathrm{add}(\mathcal{N}) = \theta_0$, $\mathrm{cov}(\mathcal{N}) = \theta_1$, $\mathfrak{b} = \mathfrak{a} = \mathrm{non}(\mathcal{M}) = \kappa$, $\mathrm{cov}(\mathcal{M}) = \mathfrak{d} = \nu$ *and* $\mathrm{non}(\mathcal{N}) = \mathfrak{c} = \lambda$.

*Proof.* The construction of the standard coherent system that forces the above is very similar to the one in the proof of Theorem 4.3. The only changes are that $S^{\mathbf{m}} := \emptyset$ and, for each $\xi \in [\lambda\rho, \lambda(\rho+1))$, $\dot{\mathbb{Q}}^{\mathbf{m}}_{\xi} := \mathbb{D}^{V_{\Delta(\xi),\xi}}$ when $\xi = \lambda\rho$, and $\dot{\mathbb{Q}}^{\mathbf{m}}_{\xi} := \dot{\mathbb{Q}}_{g(\varepsilon)}$ when $\xi = \lambda\rho + 1 + \varepsilon$ for some $\varepsilon < \lambda$. $\qquad\square$

**Theorem 4.5.** *Assume that $\theta_0 \leq \kappa \leq \mu$ are uncountable regular cardinals, $\nu \leq \lambda$ are cardinals such that $\mu \leq \nu$, $\nu^{<\kappa} = \nu$ and $\lambda^{<\theta_0} = \lambda$. Then, for each of the statements below, there is a ccc poset forcing it.*

(a) $\mathrm{MA}_{<\theta_0}$, $\mathrm{add}(\mathcal{N}) = \theta_0$, $\mathfrak{b} = \mathfrak{a} = \kappa$, $\mathrm{cov}(\mathcal{I}) = \mathrm{non}(\mathcal{I}) = \mu$ *for* $\mathcal{I} \in \{\mathcal{M}, \mathcal{N}\}$, $\mathfrak{d} = \nu$ *and* $\mathrm{cof}(\mathcal{N}) = \mathfrak{c} = \lambda$.

(b) $\mathrm{MA}_{<\theta_0}$, $\mathrm{add}(\mathcal{N}) = \theta_0$, $\mathrm{cov}(\mathcal{N}) = \kappa$, $\mathrm{add}(\mathcal{M}) = \mathrm{cof}(\mathcal{M}) = \mu$, $\mathrm{non}(\mathcal{N}) = \nu$ *and* $\mathrm{cof}(\mathcal{N}) = \mathfrak{c} = \lambda$.

(c) $\mathrm{MA}_{<\theta_0}$, $\mathrm{add}(\mathcal{N}) = \theta_0$, $\mathrm{cov}(\mathcal{N}) = \mathfrak{b} = \mathfrak{a} = \kappa$, $\mathrm{non}(\mathcal{M}) = \mathrm{cov}(\mathcal{M}) = \mu$, $\mathfrak{d} = \mathrm{non}(\mathcal{N}) = \nu$ *and* $\mathrm{cof}(\mathcal{N}) = \mathfrak{c} = \lambda$.

(d) $\mathrm{MA}_{<\theta_0}$, $\mathrm{add}(\mathcal{N}) = \theta_0$, $\mathrm{cov}(\mathcal{N}) = \mathfrak{b} = \mathfrak{a} = \mathrm{non}(\mathcal{M}) = \kappa$, $\mathrm{cov}(\mathcal{M}) = \mathfrak{d} = \mathrm{non}(\mathcal{N}) = \nu$ *and* $\mathrm{cof}(\mathcal{N}) = \mathfrak{c} = \lambda$.

**Theorem 4.6.** *Assume that $\kappa \leq \mu$ are uncountable regular cardinals, $\nu \leq \lambda$ are cardinals such that $\mu \leq \nu$, $\nu^{<\kappa} = \nu$ and $\lambda^{\aleph_0} = \lambda$. Then, for each of the statements below, there is a ccc poset forcing it.*

(a) $\mathrm{add}(\mathcal{N}) = \mathrm{cov}(\mathcal{N}) = \mathfrak{b} = \mathfrak{a} = \kappa$, $\mathrm{non}(\mathcal{M}) = \mathrm{cov}(\mathcal{M}) = \mu$, $\mathfrak{d} = \mathrm{non}(\mathcal{N}) = \mathrm{cof}(\mathcal{N}) = \nu$ *and* $\mathfrak{c} = \lambda$.

(b) $\mathrm{add}(\mathcal{N}) = \mathfrak{b} = \mathfrak{a} = \kappa$, $\mathrm{cov}(\mathcal{I}) = \mathrm{non}(\mathcal{I}) = \mu$ *for* $\mathcal{I} \in \{\mathcal{M}, \mathcal{N}\}$, $\mathfrak{d} = \mathrm{cof}(\mathcal{N}) = \nu$ *and* $\mathfrak{c} = \lambda$.

(c) $\mathrm{add}(\mathcal{N}) = \mathrm{cov}(\mathcal{N}) = \kappa$, $\mathrm{add}(\mathcal{M}) = \mathrm{cof}(\mathcal{M}) = \mu$, $\mathrm{non}(\mathcal{N}) = \mathrm{cof}(\mathcal{N}) = \nu$ *and* $\mathfrak{c} = \lambda$.

(d) $\mathrm{add}(\mathcal{N}) = \kappa$, $\mathrm{cov}(\mathcal{N}) = \mathrm{add}(\mathcal{M}) = \mathrm{cof}(\mathcal{M}) = \mathrm{non}(\mathcal{N}) = \mu$, $\mathrm{cof}(\mathcal{N}) = \nu$ *and* $\mathfrak{c} = \lambda$.

(e) $\mathrm{add}(\mathcal{N}) = \mathrm{non}(\mathcal{M}) = \mathfrak{a} = \kappa$, $\mathrm{cov}(\mathcal{M}) = \mathrm{cof}(\mathcal{N}) = \nu$ *and* $\mathfrak{c} = \lambda$.

*Moreover, if $\lambda^{<\kappa} = \lambda$, $\mathrm{MA}_{<\kappa}$ can be forced additionally at each of the items above.*

**Remark 4.7.** This method can be used to force values (even singular) to other cardinal invariants different than those from Cichoń's diagram. For instance, the results in [Mej17, Sect. 3] can be adapted to the present approach.

## 5. Bonus track: linkedness properties

The notions of $\sigma$-linked, $\sigma$-centered, $\sigma$-Fr-linked, etc., can be put into the following general framework.

**Definition 5.1.** Say that $\Gamma$ is a *linkedness property (for subsets of posets)* if $\Gamma$ is a class-function with domain the class of posets such that, for any poset $\mathbb{P}$, $\Gamma(\mathbb{P}) \subseteq \mathcal{P}(\mathbb{P})$ [3]. We define the following notions for a linkedness property $\Gamma$.

(1) $\Gamma$ is *basic* if $[\mathbb{P}]^{\leq 1} \subseteq \Gamma(\mathbb{P})$ for any poset $\mathbb{P}$.

(2) $\Gamma$ is *conic* if, for any poset $\mathbb{P}$, $P \subseteq \mathbb{P}$ and $Q \in \Gamma(\mathbb{P})$, if $P \subseteq \{p \in \mathbb{P} : \exists q \in Q(q \leq p)\}$ and $Q = \{q \in Q : \exists p \in P(q \leq p)\}$ then $P \in \Gamma(\mathbb{P})$.

(3) $\Gamma$ is a *downwards forcing linkedness property* if, for any complete embedding $\iota : \mathbb{P} \to \mathbb{Q}$ between posets, if $P \subseteq \mathbb{P}$ and $\iota[P] \in \Gamma(\mathbb{Q})$ then $P \in \Gamma(\mathbb{P})$.

(4) $\Gamma$ is an *upwards forcing linkedness property* if, for any complete embedding $\iota : \mathbb{P} \to \mathbb{Q}$ between posets, if $P \in \Gamma(\mathbb{P})$ and $\iota{\restriction}P$ is 1-1 then $\iota[P] \in \Gamma(\mathbb{Q})$.

(5) A *forcing linkedness property* is a downwards and upwards linkedness forcing property.

(6) $\Gamma$ is *appropriate* if it is a basic conic forcing linkedness property.

---

[3]Concretely, $\Gamma$ is a formula $\varphi(x, y)$ (with fixed parameters) in the language of ZF and $\Gamma(\mathbb{P}) := \{Q \subseteq \mathbb{P} : \varphi(Q, \mathbb{P})\}$.

(7) $\Gamma$ is *closed* if, for any poset $\mathbb{P}$ and $Q \subseteq Q' \subseteq \mathbb{P}$, if $Q' \in \Gamma(\mathbb{P})$ then $Q \in \Gamma(\mathbb{P})$.

(8) Let $\mu$ be an infinite cardinal. A poset $\mathbb{P}$ is $\mu$-$\Gamma$-*covered* if it can be covered by $\leq \mu$-many sets from $\Gamma(\mathbb{P})$. When $\mu = \aleph_0$ we just say $\sigma$-$\Gamma$-*covered*.

(9) Let $\theta$ be an infinite cardinal. A poset $\mathbb{P}$ is $\theta$-$\Gamma$-*Knaster* if, for any $P \subseteq \mathbb{P}$ of size $\theta$, there is some $Q \subseteq P$ of size $\theta$ such that $Q \in \Gamma(\mathbb{P})$. For $\theta = \aleph_1$ we just say $\Gamma$-*Knaster*.

(10) If $\Lambda$ is another linkedness property, say that $\Lambda$ *is stronger than* $\Gamma$ (or $\Gamma$ *is weaker than* $\Lambda$), denoted by $\Lambda \Rightarrow \Gamma$, if $\Lambda(\mathbb{P}) \subseteq \Gamma(\mathbb{P})$ for any poset $\mathbb{P}$. We say that both properties are *equivalent*, denote by $\Lambda \Leftrightarrow \Gamma$, when one is weaker and stronger than the other.

**Remark 5.2.** Let $\Gamma$ be a linkedness property and $\mu$ an infinite cardinal.

(1) If $\Gamma$ is closed, then any $\mu$-$\Gamma$-covered poset is $\mu^+$-$\Gamma$-Knaster.

(2) If $\Gamma$ is a closed conic forcing linkedness property then '$\mu$-$\Gamma$-covered' is a property of forcing notions.

(3) If $\Gamma$ is appropriate and $\mu$ is regular, then '$\mu$-$\Gamma$-Knaster' is a property of forcing notions.

(4) If a property $\Lambda$ is stronger than $\Gamma$ then any $\mu$-$\Lambda$-covered (Knaster) poset is $\mu$-$\Gamma$-covered (Knaster).

To see that closed is necessary in (1) and (2), consider the property $\Gamma_0$ defined by $Q \in \Gamma_0(\mathbb{P})$ iff $Q$ is not an antichain in $\mathbb{P}$ of size $\geq 2$. Note that $\Gamma_0$ is an appropriate linkedness property, but it is not closed. A poset $\mathbb{P}$ is $\mu$-$\Gamma_0$-covered iff either it is an antichain in itself of size $\leq \mu$, or it is not an antichain in itself. Hence, though $\omega_1^{\leq 1}$ and $\omega_1^1$ (as posets of sequences in $\omega_1$ ordered by $\supseteq$) are forcing equivalent, the first is $\sigma$-$\Gamma_0$-covered while the second is not. On the other hand, any poset is $\mu$-$\Gamma_0$-Knaster iff it is $\mu$-cc.

**Example 5.3.** The following are appropriate closed linkedness properties.

$\Gamma_{\eta\text{-cc}}$ (when $2 \leq \eta \leq \omega$): $\eta$-cc, that is, $Q \in \Gamma_{\eta\text{-cc}}(\mathbb{P})$ iff $Q$ does not contain antichains in $\mathbb{P}$ of size $\eta$.

$\Gamma_{\text{bd-cc}}$: $n$-cc for some $2 \leq n < \omega$.

$\Lambda_n$ (when $2 \leq n < \omega$): $n$-linked.

$\Lambda_\omega$: centered.

$\Gamma_{\text{cone}}$: Say that $Q \in \Gamma_{\text{cone}}(\mathbb{P})$ if there is some $q \in \mathbb{P}$ such that $\forall p \in Q(q \leq^* p)$ ([4]).

$\Lambda_{\text{Fr}}$: Frechet-linked.

It is clear that $\Gamma_{\text{cone}} \Rightarrow \Lambda_\omega \Rightarrow \Lambda_{n+1} \Rightarrow \Lambda_n \Rightarrow \Lambda_2 \Rightarrow \Gamma_{n\text{-cc}} \Rightarrow \Gamma_{n+1\text{-cc}} \Rightarrow \Gamma_{\text{bd-cc}} \Rightarrow \Gamma_{\omega\text{-cc}}$ for $2 \leq n < \omega$ (actually $\Gamma_{2\text{-cc}} \Leftrightarrow \Lambda_2$). Also $\Gamma_{\text{cone}} \Rightarrow \Lambda_{\text{Fr}} \Rightarrow \Gamma_{\omega\text{-cc}}$ and $\Lambda_2 \Rightarrow \Gamma_0$. These properties determine some well-known forcing properties, for example, '$\mu$-$\Lambda_\omega$-covered' means '$\mu$-centered', '$\mu$-$\Lambda_\omega$-Knaster' means 'precaliber $\mu$', '$\mu$-$\Lambda_2$-covered' means '$\mu$-linked', '$\mu$-$\Lambda_2$-Kanster' is the typical $\mu$-Knaster property, and '$\mu$-$\Lambda_{\text{Fr}}$-covered' is what we defined as $\mu$-Fr-linked in Definition 3.24.

By an argument similar to [HT48, Thm. 2.4] it can be proved that, if $\theta$ is regular, then any $Q \in \Gamma_{\omega\text{-cc}}(\mathbb{P})$ of size $\theta$ contains a 2-linked subset of the same size, thus $\theta$-$\Lambda_2$-Knaster is equivalent to $\theta$-$\Gamma_{\omega\text{-cc}}$-Knaster. On the other hand, Todorčević [Tod91, Tod86] constructed a $\Lambda_\omega$-Knaster (i.e. $\aleph_1$-precaliber) poset that is not $\sigma$-$\Gamma_{\omega\text{-cc}}$-covered (i.e. $\sigma$-finite-cc) and, under $\mathfrak{b} = \aleph_1$, a $\sigma$-$\Lambda_n$-covered poset that is not $\Lambda_{n+1}$-Knaster. Todorčević [Tod14] and Thümmel [Thü14] constructed $\sigma$-finite-cc posets that are not $\sigma$-$\Gamma_{\text{bd-cc}}$-covered (i.e. $\sigma$-bounded-cc).

It is clear that any Boolean algebra that admits a strictly positive fam (finitely additive measure) is $\sigma$-bounded-cc and, by Lemma 3.29, any complete Boolean algebra that admits a strictly positive $\sigma$-additive measure is $\sigma$-$\Lambda_{\mathrm{Fr}}$-covered. Note that $\mathbb{D}$ is a $\sigma$-centered poset that admits a strictly positive fam but it is not $\sigma$-Frechet-linked (otherwise it would contradict Theorem 3.30) and $\mathbb{B}_{\mathfrak{c}^+}$ is a complete Boolean algebra that admits a strictly positive $\sigma$-additive measure but it is not $\sigma$-linked (see Dow and Steprans [DS94]).

---

[4]Here, $\leq^*$ denotes the separable order of $\mathbb{P}$, that is, $q \leq^* p$ iff any condition compatible with $q$ in $\mathbb{P}$ is compatible with $p$.

Note that any poset is $\mu$-$\Gamma_{\mathrm{cone}}$-covered iff it is forcing equivalent to a poset of size $\leq \mu$, and the notion $\theta$-$\Gamma_{\mathrm{cone}}$-Knaster is equivalent to $(\theta, \theta)$-caliber. The property $\Gamma_{\mathrm{cone}}$ is the strongest of all the appropriate linkedness properties with respect to the class of separative posets, so it is morally the strongest appropriate linkedness property.

With the exception of $\Lambda_{\mathrm{Fr}}$, all the other properties (including $\Gamma_0$) are absolute for transitive models of ZFC. Recall from [Paw92] that there is a $\omega^\omega$-bounding proper poset that forces that $\mathbb{B}^V$ (random forcing from the ground model) adds a dominating real, so this poset, though $\sigma$-Frechet-linked in the ground model, is not forced to be so.

**Example 5.4.** In the work in progress [BCM] we discuss properties stronger than $\Lambda_{\mathrm{Fr}}$. Given a free filter $F$ of subsets of $\omega$, define $\Lambda_F$ such that, for any poset $\mathbb{P}$, $Q \in \Lambda_F(\mathbb{P})$ iff for any sequence $\langle p_n : n < \omega \rangle$ in $Q$ there is some $q \in \mathbb{P}$ that forces $\{n < \omega : p_n \in \dot{G}\} \in F^+$ (that is, it intersects every member of $F$), which is an appropriate closed linkedness property. Also define $\Lambda_{\mathrm{uf}}(\mathbb{P}) := \bigcap \{\Lambda_F(\mathbb{P}) : F \text{ free filter}\}$. It is clear that $\Lambda_{\mathrm{uf}} \Rightarrow \Lambda_{F'} \Rightarrow \Lambda_F \Rightarrow \Lambda_{\mathrm{Fr}}$ whenever $F \subseteq F'$. Even more, we have the following equivalences.

**Lemma 5.5.** *(a) For any free filter $F$ in $\omega$ generated by $< \mathfrak{p}$-many sets, $\Lambda_F \Leftrightarrow \Lambda_{\mathrm{Fr}}$.*
*(b) For any $\mathfrak{p}$-cc poset $\mathbb{P}$, $\Lambda_{\mathrm{uf}}(\mathbb{P}) = \Lambda_{\mathrm{Fr}}(\mathbb{P})$.*

*Proof.* Both items can be proved simultaneously. Let $\mathbb{P}$ a poset, $F$ a free filter on $\omega$ and assume that either $F$ is generated by $< \mathfrak{p}$-many sets or $\mathbb{P}$ is $\mathfrak{p}$-cc. It is enough to show that $\Lambda_{\mathrm{Fr}}(\mathbb{P}) \subseteq \Lambda_F(\mathbb{P})$. Assume that $Q \subseteq \mathbb{P}$ is Fr-linked but not in $\Lambda_F(\mathbb{P})$, so there are a countable sequence $\langle p_n : n < \omega \rangle$ in $Q$, a maximal antichain $A \subseteq \mathbb{P}$ and a sequence $\langle a_r : r \in A \rangle$ in $F$ such that each $r \in A$ is incompatible with $p_n$ for every $n \in a_r$. In any of the two cases of the hypothesis, it can be concluded that there is some pseudo-intersection $a \in [\omega]^{\aleph_0}$ of $\langle a_r : r \in A \rangle$. Hence each $r \in A$ forces $p_n \in \dot{G}$ for only finitely many $n \in a$, which means that $\mathbb{P}$ forces the same. However, as $Q$ is Fr-linked, there is some $q \in \mathbb{P}$ that forces $\exists^\infty n \in a(p_n \in \dot{G})$, a contradiction. $\square$

As a consequence of the previous result and Lemma 3.29, $\mathbb{E}$ and any complete Boolean algebra that admits a strictly positive $\sigma$-additive measure are $\sigma$-$\Lambda_{\mathrm{uf}}$-covered. In [BCM] we show that $\Lambda_{\mathrm{Fr}}$-Knaster posets do not add dominating reals. Hence, $\mathbb{D}$ becomes an example of a $\sigma$-centered poset that is not $\Lambda_{\mathrm{Fr}}$-Knaster. Also note that these properties associated with filters are not absolute.

To finish, in the general context of Definition 5.1, we provide simple conditions to understand when the FS iteration of $\theta$-$\Gamma$-Knaster (or covered) posets is $\theta$-$\Gamma$-Knaster (or covered), likewise for FS products. These conditions are summarized in the following definition, and they just represent facts extracted from the typical proofs of the iteration results for $\sigma$-linked and $\Lambda_2$-Knaster. At the end of this section, we relate the linkedness properties presented so far with the notions below.

**Definition 5.6.** Let $\Gamma$ be a linkedness property.

(1) $\Gamma$ is *productive* if, for any posets $\mathbb{P}$ and $\mathbb{Q}$, and for any $Q \subseteq \mathbb{P} \times \mathbb{Q}$, if $\mathrm{dom} Q \in \Gamma(\mathbb{P})$ and $\mathrm{ran} Q \in \Gamma(\mathbb{Q})$ then $Q \in \Gamma(\mathbb{P} \times \mathbb{Q})$.

(2) $\Gamma$ is *FS-productive* if for any sequence $\langle \mathbb{P}_i : i \in I \rangle$ of posets, $n < \omega$ and any $Q \subseteq \{p \in \prod_{i \in I}^{<\omega} \mathbb{P}_i : |\mathrm{dom}\, p| \leq n\}$ (FS product), if $\{\mathrm{dom}\, p : p \in Q\}$ forms a $\Delta$-system with root $s$ and $\{p{\upharpoonright}s : p \in Q\} \in \Gamma(\prod_{i \in s} \mathbb{P}_i)$ then $Q \in \Gamma(\prod_{i \in I}^{<\omega} \mathbb{P}_i)$.

(3) $\Gamma$ is *strongly productive* if, for any sequence $\langle \mathbb{P}_i : i \in I \rangle$ of posets, $Q \in \Gamma(\prod_{i \in I}^{<\omega} \mathbb{P}_i)$ whenever
    (i) there is some $n < \omega$ such that $Q \subseteq \{p \in \prod_{i \in I}^{<\omega} \mathbb{P}_i : |\mathrm{dom}\, p| \leq n\}$ and
    (ii) $\{p(i) : p \in Q\} \in \Gamma(\mathbb{P}_i)$ for any $i \in I$.

(4) $\Gamma$ is *two-step iterative* if, for any poset $\mathbb{P}$, any $\mathbb{P}$-name $\dot{Q}$ of a poset and any $Q \subseteq \mathbb{P} * \dot{Q}$, if $\mathrm{dom}\, Q \in \Gamma(\mathbb{P})$ and $\mathbb{P}$ forces that $\{\dot{q} : \exists p \in \dot{G}((p, \dot{q}) \in Q)\} \in \Gamma(\dot{Q})$, then $Q \in \Gamma(\mathbb{P} * \dot{Q})$.

(5) $\Gamma$ is *direct-limit iterative* if whenever
    (i) $\theta$ is an uncountable regular cardinal,
    (ii) $\langle \mathbb{P}_\alpha : \alpha \leq \delta \rangle$ is an increasing $\lessdot$-sequence of posets such that $\mathrm{cf}(\delta) = \theta$ and $\mathbb{P}_\gamma = \mathrm{limdir}_{\alpha < \gamma} \mathbb{P}_\alpha$ for any limit $\gamma \leq \delta$,
    (iii) $f : \theta \to \delta$ is increasing,
    (iv) $Q = \{p_\xi : \xi < \theta\} \subseteq \mathbb{P}_\delta$ such that each $p_\xi \in \mathbb{P}_{f(\xi+1)}$, and
    (v) for each $\xi < \theta$, $r_\xi \in \mathbb{P}_{f(\xi)}$ is a reduction of $p_\xi$,
if there is some $\gamma < \delta$ such that $r_\xi \in \mathbb{P}_\gamma$ for every $\xi < \theta$ and $\{r_\xi : \xi < \theta\} \in \Gamma(\mathbb{P}_\gamma)$, then $Q \in \Gamma(\mathbb{P}_\delta)$.

(6) $\Gamma$ is *strongly iterative* if, for any FS iteration $\mathbb{P}_\delta = \langle \mathbb{P}_\alpha, \dot{Q}_\alpha : \alpha < \delta \rangle$, $Q \in \Gamma(\mathbb{P}_\delta)$ whenever
    (i) there is some $n < \omega$ such that $Q \subseteq \{p \in \mathbb{P}_\delta : |\mathrm{dom}\, p| \leq n\}$ and
    (ii) for any $\alpha < \delta$, if $Q{\upharpoonright}\alpha \in \Gamma(\mathbb{P}_\alpha)$ then

$$\Vdash_{\mathbb{P}_\alpha} \{p(\alpha) : p{\upharpoonright}(\alpha+1) \in Q{\upharpoonright}(\alpha+1), \ p{\upharpoonright}\alpha \in \dot{G}_\alpha\} \in \Gamma(\dot{Q}_\alpha).$$

Note that any strongly productive linkedness property is both productive and FS-productive. On the other hand, if $\Gamma$ is strongly iterative, $\mathbb{P}_\delta = \langle \mathbb{P}_\alpha, \dot{Q}_\alpha : \alpha < \delta \rangle$ is a FS iteration and $Q \subseteq \{p \in \mathbb{P}_\delta : |\mathrm{dom}\, p| \leq n\}$ satisfies (6)(i),(ii) then $Q{\upharpoonright}\alpha \in \Gamma(\mathbb{P}_\alpha)$ for any $\alpha \leq \delta$. It is clear that any strongly iterative property is two-step iterative and satisfies a weak form of direct-limit iterative (which we leave implicit in the proof of Corollary 5.10).

The following is a general result about FS products.

**Theorem 5.7.** *Let $\mu$ be an infinite cardinal, $\theta$ an uncountable regular cardinal, and let $\Gamma$ be an appropriate linkedness property.*

(a) *If $\Gamma$ is productive then any finite product of $\mu$-$\Gamma$-covered sets is $\mu$-$\Gamma$-covered.*

(b) *If $\Gamma$ is closed and productive then and any finite product of $\theta$-$\Gamma$-Knaster posets is $\theta$-$\Gamma$-Knaster.*

(c) *If $\Gamma$ is FS-productive and $\langle \mathbb{P}_i : i \in I \rangle$ is a sequence of $\theta$-$\Gamma$-Knaster posets, then $\prod_{i \in I}^{<\omega} \mathbb{P}_i$ is $\theta$-$\Gamma$-Knaster iff $\prod_{i \in s} \mathbb{P}_i$ is $\theta$-$\Gamma$-Knaster for every $s \in [I]^{<\aleph_0}$.*

(d) *If $\Gamma$ is strongly productive, $\langle \mathbb{P}_i : i \in I \rangle$ is a sequence of $\mu$-$\Gamma$-covered posets and $|I| \leq 2^\mu$, then $\prod_{i \in I}^{<\omega} \mathbb{P}_i$ is $\mu$-$\Gamma$-covered.*

*Proof.* Items (a),(b) are easy and (c) follows by a classical $\Delta$-system argument. Item (d) uses the following result.

**Lemma 5.8** (Engelking and Karłowicz [EK65])**.** *If $\mu$ is an infinite cardinal and $I$ is a set of size $\leq 2^\mu$ then there exists a set $H \subseteq \mu^I$ of size $\leq \mu$ such that any finite partial function from $I$ to $\mu$ is extended by some member of $H$.*

For each $i \in I$ choose a sequence $\langle Q_{i,\zeta} : \zeta < \mu \rangle$ of non-empty sets in $\Gamma(\mathbb{P}_i)$ that covers $\mathbb{P}_i$. Let $H$ be as in Lemma 5.8. By Definition 5.6(3), the set $Q_{h,n}^* := \{p \in \prod_{i \in I}^{<\omega} Q_{i,h(i)} : |\mathrm{dom}\, p| = n\}$ is in $\Gamma\left(\prod_{i \in I}^{<\omega} \mathbb{P}_i\right)$ and it is clear that $\langle Q_{h,n}^* : h \in H,\ n < \omega \rangle$ covers $\prod_{i \in I}^{<\omega} \mathbb{P}_i$. $\square$

Now we turn to a general result about FS iterations.

**Theorem 5.9.** *Let $\mu$ be an infinite cardinal, $\theta$ an uncountable regular cardinal, and let $\Gamma$ be an appropriate linkedness property.*

(a) *If $\Gamma$ is two-step iterative, $\mathbb{P}$ is $\mu$-$\Gamma$-covered, and $\mathbb{P}$ forces that $\Gamma$ is basic and that $\dot{\mathbb{Q}}$ is a $|\mu|$-$\Gamma$-covered poset, then $\mathbb{P} * \dot{\mathbb{Q}}$ is $\mu$-$\Gamma$-covered.*

(b) *Let $\mathbb{P}$ be a $\theta$-$\Gamma$-Knaster poset and let $\dot{\mathbb{Q}}$ be a $\mathbb{P}$-name of a $\theta$-$\Gamma$-Knaster poset. Assume in addition that either $\Gamma = \Gamma_0$, or $\Gamma$ is closed, two-step iterative, $\mathbb{P}$ is $\theta$-cc and $\mathbb{P}$ forces that $\Gamma$ is closed and basic. Then $\mathbb{P} * \dot{\mathbb{Q}}$ is $\theta$-$\Gamma$-Knaster.*

(c) *If $\Gamma$ is direct-limit iterative, $\delta$ is a limit ordinal and $\langle \mathbb{P}_\alpha : \alpha \leq \delta \rangle$ is an increasing $\lessdot$-sequence of $\theta$-$\Gamma$-Knaster posets such that $\mathbb{P}_\gamma = \mathrm{limdir}_{\alpha < \gamma} \mathbb{P}_\alpha$ for any limit $\gamma \leq \delta$, then $\mathbb{P}_\delta$ is $\theta$-$\Gamma$-Knaster.*

(d) *If $\Gamma$ is strongly iterative and any poset forces that $\Gamma$ is still basic, then any FS iteration of length $< (2^\mu)^+$ of $\mu$-$\Gamma$-covered posets is $\mu$-$\Gamma$-covered.*

*Proof.* To see (a), it is enough to show that, for any poset $\mathbb{P}$ that forces $\Gamma$ to be still basic and any $\mathbb{P}$-name $\dot{\mathbb{Q}}$ for a poset, if $P \in \Gamma(\mathbb{P})$ and $\dot{Q}$ is a $\mathbb{P}$-name of a non-empty set in $\Gamma(\dot{\mathbb{Q}})$, then $P * \dot{Q} := \{(p, \dot{q}) \in \mathbb{P} * \dot{\mathbb{Q}} : p \in P \text{ and } p \Vdash \dot{q} \in \dot{Q}\}$ is in $\Gamma(\mathbb{P} * \dot{\mathbb{Q}})$. Clearly $\mathrm{dom}(P * \dot{Q}) = P$. On the other hand, any $p_0 \in P$ forces that $\dot{R} := \{\dot{q} : \exists p \in \dot{G}((p, \dot{q}) \in P * \dot{Q})\} = \dot{Q}$, so $\mathbb{P}$ forces that $\dot{R}$ is either $\dot{Q}$ or the empty set, so $\dot{R} \in \Gamma(\dot{\mathbb{Q}})$. By Definition 5.6(4), it follows that $P * \dot{Q} \in \Gamma(\mathbb{P} * \dot{\mathbb{Q}})$.

Item (b) is well-known when $\Gamma = \Gamma_0$, so assume that $\Gamma$ is closed, two-step iterative, $\mathbb{P}$ is $\theta$-cc and $\mathbb{P}$ forces that $\Gamma$ is closed and basic. Let $\{(p_\alpha, \dot{q}_\alpha) : \alpha < \theta\} \subseteq \mathbb{P} * \dot{\mathbb{Q}}$. As $\mathbb{P}$ forces $\dot{\mathbb{Q}}$ to be $\theta$-$\Gamma$-Knaster, there is some $\mathbb{P}$-name $\dot{K}$ for a subset of $\theta$ such that $\mathbb{P}$ forces that " whenever $|\{\alpha < \theta : p_\alpha \in \dot{G}\}| = \theta$, $\dot{K} \subseteq \{\alpha < \theta : p_\alpha \in \dot{G}\}$ has size $\theta$ and $\{\dot{q}_\alpha : \alpha \in \dot{K}\} \in \Gamma(\dot{\mathbb{Q}})$, otherwise $\dot{K} = \emptyset$". Set $K_0 := \{\alpha < \theta : \nVdash \alpha \notin \dot{K}\}$, which has size $\theta$ (otherwise $\mathbb{P}$ would force that $|\{\alpha < \theta : p_\alpha \in \dot{G}\}| < \theta$, which contradicts that $\mathbb{P}$ is $\theta$-cc). For each $\alpha \in K_0$ choose an $r_\alpha \leq p_\alpha$ that forces $\alpha \in \dot{K}$. Hence, there is some $K_1 \subseteq K_0$ of size $\theta$ such that $\{r_\alpha : \alpha \in K_1\} \in \Gamma(\mathbb{P})$.

As $\Gamma$ is conic and $r_\alpha \leq p_\alpha$ for any $\alpha \in K_1$, it is enough to show that $Q := \{(r_\alpha, \dot{q}_\alpha) : \alpha \in K_1\} \in \Gamma(\mathbb{P} * \dot{\mathbb{Q}})$. It is clear that $\mathrm{dom}\, Q \in \Gamma(\mathbb{P})$. On the other hand, $\mathbb{P}$ forces that $\dot{R} := \{\dot{q}_\alpha : r_\alpha \in \dot{G},\ \alpha \in K_1\} \subseteq \{\dot{q}_\alpha : \alpha \in \dot{K}\}$, so $\dot{R} \in \Gamma(\dot{\mathbb{Q}})$ because $\Gamma$ is closed. As $\Gamma$ is two-step iterative, we are done.

Now we show (c). Let $\{p_\xi : \xi < \theta\} \subseteq \mathbb{P}_\delta$. If $\mathrm{cf}(\delta) \neq \theta$ then there are some $\alpha < \delta$ and a $K \subseteq \theta$ of size $\theta$ such that $\{p_\xi : \xi \in K\} \subseteq \mathbb{P}_\alpha$, so there is some $K' \subseteq K$ of size $\theta$ such that $\{p_\xi : \xi \in K\} \in \Gamma(\mathbb{P}_\alpha)$ (note that $\Gamma(\mathbb{P}_\alpha) \subseteq \Gamma(\mathbb{P}_\delta)$ because $\Gamma$ is appropriate). Assume that $\mathrm{cf}(\delta) = \theta$ and choose an increasing continuous cofinal function $g : \theta \to \delta$ such that each $g(\xi)$ is a limit ordinal. For each $\xi < \theta$ choose a reduction $r_\xi \in \mathbb{P}_{g(\xi)}$ of $p_\xi$. As $g(\xi)$ is limit, there is some $h(\xi) < \xi$ such that $r_\xi \in \mathbb{P}_{g(h(\xi))}$. Hence, by Fodor's Lemma, there is some stationary set $S \subseteq \theta$ such that $h[S] = \{\eta\}$ for some $\eta < \theta$, that is, $r_\xi \in \mathbb{P}_{g(\eta)}$ for every $\xi \in S$. By recursion define $j : \theta \to S$ increasing such that $j(0) > \eta$ and, for any $\zeta < \theta$, $p_{j(\zeta)} \in \mathbb{P}_{g(j(\zeta+1))}$. As $\mathbb{P}_{g(\eta)}$ is $\theta$-$\Gamma$-Knaster, there is some $K \subseteq \theta$ of size $\theta$ such that $\{r_{j(\zeta)} : \zeta \in K\} \in \Gamma(\mathbb{P}_{g(\eta)})$. Let $i : \theta \to K$ be the increasing enumeration of $K$.

Put $f := g \circ j \circ i$ and $\gamma := g(\eta)$. Note that $\{p_{j(i(\beta))} : \beta < \theta\}$, $\{r_{j(i(\beta))} : \beta < \theta\}$, $f$ and $\gamma$ satisfy the conditions of Definition 5.6(5) so, as $\Gamma$ is direct-limit iterative, $\{p_{j(i(\beta))} : \beta < \theta\} \in \Gamma(\mathbb{P}_\delta)$.

To finish, we show (d). Let $\delta < (2^\mu)^+$ and let $\mathbb{P}_\delta = \langle \mathbb{P}_\alpha, \dot{\mathbb{Q}}_\alpha : \alpha < \delta \rangle$ be a FS iteration of $\mu$-$\Gamma$-covered sets. For each $\alpha < \delta$ choose a sequence $\langle \dot{Q}_{\alpha,\zeta} : \zeta < \mu \rangle$ of $\mathbb{P}_\alpha$-names of sets in $\Gamma(\dot{\mathbb{Q}}_\alpha)$ that is forced to cover $\dot{\mathbb{Q}}_\alpha$. For $\alpha \leq \delta$ define $\mathbb{P}_\alpha^* \subseteq \mathbb{P}_\alpha$ such that $p \in \mathbb{P}_\alpha^*$ iff $p \in \mathbb{P}_\alpha$ and, for any $\xi \in \operatorname{dom} p$, there is some $\zeta < \mu$ such that $p{\restriction}\xi \Vdash_{\mathbb{P}_\xi} p(\xi) \in \dot{Q}_{\xi,\zeta}$. By induction it can be proved that $\mathbb{P}_\alpha^*$ is dense in $\mathbb{P}_\alpha$.

Now choose $H$ as in Lemma 5.8 and, for each $h \in H$ and $n < \omega$, define $Q_{h,n}$ as the set of $p \in \mathbb{P}_\delta^*$ such that $|\operatorname{dom} p| \leq n$ and, for any $\alpha \in \operatorname{dom} p$, $p{\restriction}\alpha \Vdash_{\mathbb{P}_\alpha} p(\alpha) \in \dot{Q}_{\alpha,h(\alpha)}$. It is clear that $\langle Q_{h,n} : h \in H, n < \omega \rangle$ covers $\mathbb{P}_\delta^*$, so it remains to show that $Q_{h,n} \in \Gamma(\mathbb{P}_\delta)$. If $\alpha < \delta$ and $Q_{h,n}{\restriction}\alpha \in \Gamma(\mathbb{P}_\alpha)$ then a similar argument as in (a) shows that $\mathbb{P}_\alpha$ forces $\{p(\alpha) : p{\restriction}(\alpha+1) \in Q_{h,n}{\restriction}(\alpha+1), \ p{\restriction}\alpha \in \dot{G}\} \in \Gamma(\dot{\mathbb{Q}}_\alpha)$. Therefore, as $\Gamma$ is strongly iterative, $Q_{h,n} \in \Gamma(\mathbb{P}_\delta)$. $\qquad\square$

**Corollary 5.10.** *Let $\theta$ an uncountable regular cardinal and assume that $\Gamma$ is either*

(i) *$\Gamma_0$ or*

(ii) *a closed appropriate linkedness property that is closed and basic in any generic extension, and that it is either strongly iterative, or two-step and direct-limit iterative.*

*Then any FS iteration of $\theta$-$\Gamma$-Knaster $\theta$-cc posets is $\theta$-$\Gamma$-Knaster.*

*Proof.* Case (i) and case (ii) when $\Gamma$ is two-step and direct-limit iterative follow directly from Theorem 5.9. Case (ii) when $\Gamma$ is strongly iterative is a bit similar but requires a bit more work. If $\langle \mathbb{P}_\alpha, \dot{\mathbb{Q}}_\alpha : \alpha < \delta \rangle$ is a FS iteration of $\theta$-$\Gamma$-Knaster $\theta$-cc posets, it is enough to show by induction on $\alpha \leq \delta$ that, for any sequence $\langle p_\beta : \beta < \theta \rangle$ in $\mathbb{P}_\alpha$ there are some $K \subseteq \theta$ of size $\theta$ and some sequence $\langle r_\beta : \beta \in K \rangle$ in $\mathbb{P}_\alpha$ that satisfies (i) and (ii) of Definition 5.6(6) (with respect to $\mathbb{P}_\alpha$) and such that $r_\beta \leq p_\beta$ for any $\beta \in K$. The successor step is exactly like the proof of Theorem 5.9(b) and the limit step is very similar to Theorem 5.9(c). We just look at the case $\operatorname{cf}(\alpha) = \theta$. Let $\langle p_\beta : \beta < \theta \rangle$ be a sequence in $\mathbb{P}_\alpha$. Exactly like in the proof of Theorem 5.9(c), we can find a $\gamma < \alpha$, a $K_0 \subseteq \theta$ of size $\theta$ and an increasing function $f : K_0 \to \alpha \smallsetminus \gamma$ such that, for each $\beta \in K_0$, $p_\beta{\restriction}f(\beta) \in \mathbb{P}_\gamma$ and $p_\beta \in \mathbb{P}_{f(\beta+1)}$. Even more, we may assume that there is some $n_1 < \omega$ such that $|\operatorname{dom} p_\beta \smallsetminus f(\beta)| = n_1$ for all $\beta \in K$. By the inductive hypothesis, there are $K \subseteq K_0$ of size $\theta$ and a sequence $\langle r_\beta^0 : \beta \in K \rangle$ of conditions in $\mathbb{P}_\gamma$ that satisfies (i) (for some $n_0 < \omega$) and (ii) of Definition 5.6(6) (with respect to $\mathbb{P}_\gamma$) and such that $r_\beta^0 \leq p_\beta{\restriction}f(\beta)$ for any $\beta \in K$. The set $\{r_\beta^0 \cup p_\beta{\restriction}(f(\beta+1) \smallsetminus f(\beta)) : \beta \in K\}$ is as required. $\qquad\square$

**Remark 5.11.** Table 1 illustrates which productive or iterative notions are satisfied by the linkedness properties discussed so far.

We explain some of the facts indicated in the table. First, we show that $\Gamma_{\omega\text{-cc}}$ is strongly productive. Let $Q \subseteq \{p \in \prod_{i \in I}^{<\omega} \mathbb{P}_i : |\operatorname{dom} p| \leq n\}$ and assume that $\{p(i) : p \in Q\} \in \Gamma_{\omega\text{-cc}}(\mathbb{P}_i)$ for every $i \in I$. Fix a countable sequence $\langle p_k : k < \omega \rangle$ in $Q$. As the size of the domains of the members of the sequence are bounded by $n$, we can find a $W \in [\omega]^{\aleph_0}$ such that $\langle p_k : k \in W \rangle$ forms a $\Delta$-system with root $R$. By Ramsey's Theorem it can be proved that $\{p_k{\restriction}R : k \in W\}$ is not an antichain, so neither is $\langle p_k : k < \omega \rangle$.

Ramsey's Theorem also implies that $\Gamma_{\text{bd-cc}}$ is productive. However, it is not FS-productive (consider the set of conditions with domain of size 1 of the FS product $\prod_{i \in \omega}^{<\omega} \mathbb{P}_i$ where each $\mathbb{P}_i$ is an antichain of size $i+1$). The following example indicates that both $\Gamma_{\omega\text{-cc}}$ and $\Gamma_{\text{bd-cc}}$ are not two-step iterative. Consider $P := \{(p_n, \dot{q}_n) : n < \omega\} \subseteq \mathbb{C} * \dot{\mathbb{C}}$ such that each $p_n$ is a sequence of zeros of length $n+1$, $p_n \Vdash \dot{q}_n = \langle n \rangle$ but $r \Vdash \dot{q}_n = \langle \ \rangle$ for any $r \in \mathbb{C}$ incompatible with $p_n$. Though $P$ is an antichain in $\mathbb{C} * \dot{\mathbb{C}}$, $\mathrm{dom}P$ is centered and $\mathbb{C}$ forces that $\{\dot{q}_n : p_n \in \dot{G}, \ n < \omega\}$ is a finite antichain.

| | $\Gamma_0$ | $\Gamma_{n\text{-cc}}$ $(3 \leq n)$ | $\Gamma_{\omega\text{-cc}}$ | $\Gamma_{\text{bd-cc}}$ | $\Lambda_n$ $(2 \leq n)$ | $\Lambda_\omega$ | $\Gamma_{\text{cone}}$ | $\Lambda_{\text{Fr}}$ | $\Lambda_F$ | $\Lambda_{\text{uf}}$ |
|---|---|---|---|---|---|---|---|---|---|---|
| Prod. | $\times$ | $\times$ | $\bigcirc$ | $\bigcirc$ | $\bigcirc$ | $\bigcirc$ | $\bigcirc$ | ? | ? | $\bigcirc$ |
| FS Prod. | $\bigcirc$ | $\bigcirc$ | $\bigcirc$ | $\bigcirc$ | $\bigcirc$ | $\bigcirc$ | $\times$ | $\bigcirc$ | $\bigcirc$ | $\bigcirc$ |
| Str. Prod. | $\times$ | $\times$ | $\bigcirc$ | $\times$ | $\bigcirc$ | $\bigcirc$ | $\times$ | ? | ? | ? |
| Two-step it. | $\bigcirc$ | $\times$ | $\times$ | $\times$ | $\bigcirc$ | $\bigcirc$ | $\bigcirc$ | $\bigcirc$ | ? | $\bigcirc$ |
| Dir.-lim. it. | $\bigcirc$ | $\bigcirc$ | $\bigcirc$ | $\bigcirc$ | $\bigcirc$ | $\bigcirc$ | $\times$ | $\bigcirc$ | $\bigcirc$ | $\bigcirc$ |
| Str. it | $\bigcirc$ | $\times$ | $\times$ | $\times$ | $\bigcirc$ | $\bigcirc$ | $\times$ | $\bigcirc$ | ? | ? |

Table 1. A circle means that the linkedness property satisfies the corresponding productive or iterative notion (see Definition 5.6) on the left, an $\times$ means that such notion is not satisfied, and a question mark means unclear.

As $\mathbb{C}_{\omega_1}$ is uncountable and it has not $(\aleph_1, \aleph_0)$-precaliber, Theorems 5.7(c),(d) and 5.9(c),(d) cannot be applied to $\Gamma_{\text{cone}}$. Hence, $\Gamma_{\text{cone}}$ does not satisfy the properties indicated with $\times$ in the table.

It is unclear whether $\Lambda_F$ is productive in general, but it is proved in [BCM] that it is when $F$ is an ultrafilter. Therefore, $\Lambda_{\text{uf}}$ is productive. By a $\Delta$-system argument, $\Lambda_{\text{Fr}}$ is strongly iterative. To see this, assume that $\mathbb{P}_\delta$, $Q$ and $n < \omega$ satisfy the conditions in Definition 5.6(6). It is enough to show, by induction on $\alpha \leq \delta$, that $Q \upharpoonright \alpha \in \Lambda_{\text{Fr}}(\mathbb{P}_\alpha)$. Since $\Lambda_{\text{Fr}}$ is two-step iterative, we only need to prove the limit step. Let $\langle p_k : k < \omega \rangle$ be a sequence in $Q \upharpoonright \alpha$. As $|\mathrm{dom}p_k| \leq n$ for any $k < \omega$, there is some infinite $W \subseteq \omega$ such that $\{\mathrm{dom}p_k : k \in W\}$ forms a $\Delta$-system with root $R$, and there is some $\xi < \alpha$ such that $R \subseteq \xi$. By the inductive hypothesis, $Q \upharpoonright \xi \in \Lambda_{\text{Fr}}(\mathbb{P}_\xi)$, so there is some $q \in \mathbb{P}_\xi$ that forces $\exists^\infty k \in W(p_k \upharpoonright \xi \in \dot{G}_\xi)$. Therefore, it can be proved that $q$ forces (in $\mathbb{P}_\alpha$) that $\exists^\infty k \in W(p_k \in \dot{G}_\alpha)$.

By a similar argument, if $\Lambda_{\text{Fr}}$ were productive then it would be strongly productive. In particular, by Lemma 5.5, $\Lambda_{\text{Fr}}$ restricted to the class of Knaster posets is strongly productive, so Theorem 5.7 is valid for $\Lambda_{\text{Fr}}$ for FS products of Knaster posets (or just FS products that have ccc).

**Acknowledgments.** This work was supported by grant no. IN201711, Dirección Operativa de Investigación, Institución Universitaria Pascual Bravo, and by the Grant-in-Aid for Early Career Scientists 18K13448, Japan Society for the Promotion of Science.

The author would like to thank Miguel Cardona for the very useful discussions that helped this work to take its final form. He is also very grateful with the anonymous referee for his/her very useful comments, specially for asking whether any $\sigma$-Fr-linked poset is **Lc**-good (which is answered negative with a counterexample in Remark 3.31).

CREATIVE SCIENCE COURSE (MATHEMATICS), FACULTY OF SCIENCE, SHIZUOKA UNIVERSITY, OHYA 836, SURUGA-KU, 422-8529 SHIZUOKA, JAPAN
    *E-mail address*: `diego.mejia@shizuoka.ac.jp`
    *URL*: `http://www.researchgate.com/profile/Diego_Mejia2`